\documentclass[12pt,reqno]{amsart}

\usepackage{geometry}
\usepackage{amsmath,amssymb,mathtools}
\usepackage{mathrsfs,eucal}
\usepackage{bm}

\usepackage{graphicx}
\usepackage{multirow,bigdelim}

\usepackage{enumitem}

\usepackage[colorlinks=true, citecolor=blue, linkcolor=blue, urlcolor=blue]{hyperref}

\allowdisplaybreaks[1]

\makeatletter
\@namedef{subjclassname@2020}{\textup{2020} Mathematics Subject Classification}
\makeatother

\providecommand{\doi}[1]{}
\DeclareUrlCommand\doifont{\urlstyle{tt}}
\renewcommand{\doi}[1]{%
  \href{https://doi.org/#1}{\doifont{https://doi.org/#1}}}

\renewcommand{\MR}[1]{%
  \href{https://mathscinet.ams.org/mathscinet-getitem?mr=#1}%
       {\texttt{MR#1}}}

\newtheorem{lemma}{Lemma}[section]
\newtheorem{theorem}[lemma]{Theorem}
\newtheorem{proposition}[lemma]{Proposition}
\newtheorem{corollary}[lemma]{Corollary}
\newtheorem{question}[lemma]{Question}
\theoremstyle{definition}
\newtheorem{definition}[lemma]{Definition}
\newtheorem{example}[lemma]{Example}
\newtheorem{remark}[lemma]{Remark}

\numberwithin{equation}{section}

\newcommand{\bdf}{\begin{definition}}
\newcommand{\edf}{\end{definition}}
\newcommand{\blem}{\begin{lemma}}
\newcommand{\elem}{\end{lemma}}
\newcommand{\bthm}{\begin{theorem}}
\newcommand{\ethm}{\end{theorem}}
\newcommand{\bpf}{\begin{proof}}
\newcommand{\epf}{\end{proof}}
\newcommand{\bprop}{\begin{proposition}}
\newcommand{\eprop}{\end{proposition}}
\newcommand{\bcor}{\begin{corollary}}
\newcommand{\ecor}{\end{corollary}}
\newcommand{\brem}{\begin{remark}}
\newcommand{\erem}{\end{remark}}
\newcommand{\bquest}{\begin{question}}
\newcommand{\equest}{\end{question}}
\newcommand{\bex}{\begin{example}}
\newcommand{\eex}{\end{example}}

\newcommand{\benu}{\begin{enumerate}\renewcommand{\labelenumi}{{\rm (\arabic{enumi})}}\renewcommand{\itemsep}{0pt}}
\newcommand{\eenu}{\end{enumerate}}

\newcommand{\C}{\mathbb{C}}

\newcommand{\bB}{\mathbb{B}}

\newcommand{\cA}{\mathcal{A}}
\newcommand{\cB}{\mathcal{B}}

\newcommand{\cD}{\mathcal{D}}

\newcommand{\cH}{\mathcal{H}}
\newcommand{\cK}{\mathcal{K}}
\newcommand{\cL}{\mathcal{L}}
\newcommand{\cM}{\mathcal{M}}
\newcommand{\cN}{\mathcal{N}}

\newcommand{\cP}{\mathcal{P}}

\newcommand{\cR}{\mathcal{R}}

\newcommand{\cU}{\mathcal{U}}
\newcommand{\cV}{\mathcal{V}}

\DeclareMathOperator{\id}{id}

\begin{document}
\title{Operator-valued maximal $f$-divergences for completely positive maps}
\author{Rui OKAYASU}
\address{Department of Mathematics Education, Osaka Kyoiku University, Kashiwara, Osaka 582-8582, JAPAN}
\email{rui@cc.osaka-kyoiku.ac.jp}
%\date{\today}
\thanks{The author is partially supported by JSPS KAKENHI Grant Number JP26K06843.}

\begin{abstract}
We introduce an operator-valued maximal $f$-divergence for completely
positive (CP) maps between von Neumann algebras, associated with an operator
convex function $f$ on $(0,+\infty)$.
The construction takes values in the extended
lower-semibounded self-adjoint part of the codomain von Neumann algebra.
We show that it can be computed from any CP map dominating both maps
up to positive constants, and we prove joint subadditivity,
monotonicity under unital precomposition and normal postcomposition,
a martingale convergence theorem for normal CP maps,
and joint lower semicontinuity of its values at each normal positive
functional with respect to the point-$\sigma$-weak topology.

For normal positive functionals, 
our construction recovers Hiai's maximal $f$-divergence.
As a consequence, Hiai's maximal $f$-divergence is jointly lower semicontinuous in the $\sigma(\cM_*,\cM)$-topology for every von Neumann algebra $\cM$, which answers a question left open by Hiai.

For $\eta(t)=t\log t$, we obtain an operator-valued
Belavkin--Staszewski (BS) relative entropy.
Moreover, for normal channels with $\sigma$-finite codomain,
we prove that the BS channel divergence of Hollands and Ranallo
coincides with the extended norm of our operator-valued divergence.

Finally, we give finite-dimensional examples and obtain an explicit formula for a finite-index
conditional expectation.
In the BS case, this formula
reduces to the logarithm of the Jones--Kosaki index.
\end{abstract}

\maketitle

%%%%%%%%%%%%%%%%%%%%%%%%%%%%%%%%%%%%%%%%%%%%%%%%%%%%%%%%%%%%%%%%%%%%%%%%%%%%%%%%
\section{Introduction}
%%%%%%%%%%%%%%%%%%%%%%%%%%%%%%%%%%%%%%%%%%%%%%%%%%%%%%%%%%%%%%%%%%%%%%%%%%%%%%%%

Quantum divergences quantify the distinguishability of quantum states and play a fundamental role in quantum information theory and operator algebras.
For positive matrices $\rho$ and $\sigma$, Matsumoto introduced the maximal $f$-divergence $\widehat S_f(\rho\|\sigma)$ for any convex function $f$ on $(0,+\infty)$ \cite{mat}.
In this paper, we restrict our attention to the subclass of operator convex functions.
For $\eta(t)=t\log t$, the maximal $\eta$-divergence reduces to the Belavkin--Staszewski (BS) relative entropy introduced in \cite{bs}.
Hiai subsequently extended maximal $f$-divergences associated with operator convex functions to normal positive functionals on general von Neumann algebras \cite{hiai_2019}.

Completely positive (CP) maps provide a natural framework for extending such divergences from states to quantum transformations.
Existing channel divergences are typically scalar-valued and are obtained by stabilizing a state divergence with ancillary systems and optimizing over input states;
see, for example,
\cite{wilde-berta-hirche-kaur_2020,fang-fawzi_2021}.
In the BS case, Hollands and Ranallo introduced a
completely stabilized channel divergence and applied it to complexity in algebraic quantum field theory \cite{h-r_2023}.
In particular, they related the divergence of the identity map from the minimal conditional expectation to the Jones--Kosaki index.

The purpose of this paper is to introduce an operator-valued maximal
$f$-divergence for CP maps between von Neumann algebras.
Let $\cM$ and $\cN$ be von Neumann algebras, with $\cN$ represented on a Hilbert space $\cH$.
Unless otherwise stated, the CP maps considered
below are not assumed to be normal.
For $\Phi,\Psi\colon\cM\to\cN$,
let $(\pi,V,\cK)$ be a minimal Stinespring representation of $\Phi+\Psi$,
and let $A',B'\in\pi(\cM)'$ be the corresponding Radon--Nikodym derivatives.
Using the Pusz--Woronowicz (PW) functional calculus of \cite{huw_2022}, we define
\[
 \widehat S_f(\Phi\|\Psi)
 \coloneqq
 V^*\phi_f(A',B')V,
\]
where $\phi_f$ is the perspective associated with $f$.

Since $\phi_f(A',B')$ may be unbounded,
the divergence takes values in the extended lower-semibounded self-adjoint part $\widehat{\cN}_{\mathrm{lb}}$ introduced in \cite{huw_2022}.
We prove that it can be computed using any CP map $\Gamma$ such that
$c_\Phi\Gamma-\Phi$ and $c_\Psi\Gamma-\Psi$ are CP for some constants
$c_\Phi,c_\Psi>0$.
Unlike stabilized channel divergences,
our construction retains operator-valued information in the codomain.

We establish several fundamental properties of the operator-valued maximal $f$-divergence,
including joint subadditivity,
the transpose relation,
monotonicity under unital precomposition and normal
postcomposition by CP maps,
and a martingale convergence theorem for
normal CP maps.
We also characterize when
$\widehat S_f(\Phi\|\Psi)$ has a dense domain and give sufficient conditions for it to be bounded.
These criteria are formulated in terms of the boundary values $f(+0)$ and $f'(+\infty)$ and the absolute
continuity relations between $\Phi$ and $\Psi$,
using the Lebesgue-type decomposition of CP maps developed in \cite{oka}.

A principal result is joint lower semicontinuity in the
point-$\sigma$-weak topology.
If $\Phi_i(a)\to\Phi(a)$ and $\Psi_i(a)\to\Psi(a)$ $\sigma$-weakly for every $a\in\cM$,
then
\[
 \widehat S_f(\Phi\|\Psi)(\psi)
 \leq
 \liminf_i
 \widehat S_f(\Phi_i\|\Psi_i)(\psi),
 \qquad
 \psi\in\cN_*^+.
\]
For normal positive functionals,
our definition agrees with Hiai's maximal $f$-divergence.
Consequently, Hiai's maximal $f$-divergence is jointly $\sigma(\cM_*,\cM)$-lower semicontinuous for every von Neumann algebra $\cM$, 
answering the question raised in \cite[Remark~5.7]{hiai_2019}.

For $\eta(t)=t\log t$,
we write
\[
 D_{\mathrm{BS}}^{\mathrm{op}}(\Phi\|\Psi)
 \coloneqq
 \widehat S_\eta(\Phi\|\Psi).
\]
Using the weighted geometric means of CP maps developed in \cite{oka},
we obtain an infinitesimal-limit representation of this
operator-valued BS relative entropy.
Moreover, for normal channels with $\sigma$-finite
codomain, we prove
\[
 D_{\mathrm{BS}}^{\mathrm{HR}}(\Phi\|\Psi)
 =
 \left\|
 D_{\mathrm{BS}}^{\mathrm{op}}(\Phi\|\Psi)
 \right\|_{\mathrm{ext}},
\]
so that the Hollands--Ranallo scalar divergence is recovered by taking the extended norm of our operator-valued divergence.

We finally give finite-dimensional examples and compute the divergence for finite-index conditional expectations.
The examples include Schur multipliers and congruence maps; the latter example shows that CP maps can be
mutually singular even when their values at the identity are invertible and arbitrarily close.
If $\cN\subseteq\cM$ is an inclusion of $\sigma$-finite factors and $E\colon\cM\to\cN$ is a faithful normal conditional expectation with finite Jones--Kosaki index
$\lambda=\operatorname{Ind}E$,
then, regarding $E$ as an $\cM$-valued CP map,
we prove
\[
 \widehat S_f(\id_{\cM}\|E)
 =
 \left\{
 \frac{f(\lambda)}{\lambda}
 +
 \left(1-\frac{1}{\lambda}\right)f(+0)
 \right\}1_{\cM}
\]
with the convention $0\cdot(+\infty)=0$.
In particular,
\[
 D_{\mathrm{BS}}^{\mathrm{op}}(\id_{\cM}\|E)
 =
 (\log\operatorname{Ind}E)1_{\cM}.
\]
For the minimal conditional expectation $E_0$,
the extended norm identity recovers the index formula of Hollands and Ranallo
\cite[Proposition~3.21]{h-r_2023}
\[
 D_{\mathrm{BS}}^{\mathrm{HR}}(\id_{\cM}\|E_0)
 =
 \log[\cM:\cN],
\]
where $[\cM:\cN]=\operatorname{Ind}E_0$ is the minimal index.

%%%%%%%%%%%%%%%%%%%%%%%%%%%%%%%%%%%%%%%%%%%%%%%%%%%%%%%%%%%%%%%%%%%%%%%%%%%%%%%%%%%%%%%%%%%%%%%%%%%%%%%%%%%%%%%%%%%%%%%%%%%%%%%%%%%%%%%%%%%%%%%%%%%%%%%%%%%%%%%%%%%%%%%%%%%%%%%%%%%%%%%%%%%%%%%%%%%%%%%%%%%%%%%%%%%%%%%%%%%%%%%%%%%%%%%%%%%%%%%%
\section{Preliminaries}
%%%%%%%%%%%%%%%%%%%%%%%%%%%%%%%%%%%%%%%%%%%%%%%%%%%%%%%%%%%%%%%%%%%%%%%%%%%%%%%%%%%%%%%%%%%%%%%%%%%%%%%%%%%%%%%%%%%%%%%%%%%%%%%%%%%%%%%%%%%%%%%%%%%%%%%%%%%%%%%%%%%%%%%%%%%%%%%%%%%%%%%%%%%%%%%%%%%%%%%%%%%%%%%%%%%%%%%%%%%%%%%%%%%%%%%%%%%%%%%%

Throughout this article, $\cH$ denotes a Hilbert space, and
$\bB(\cH)$ denotes the $C^*$-algebra of all bounded linear
operators on $\cH$.
Let $\cN$ be a von Neumann algebra acting on $\cH$.
We write $\cN_{\mathrm{sa}}$, $\cN_+$, $\cN_{++}$, and $\cU(\cN)$ for the sets
of self-adjoint, positive, positive invertible elements, and unitaries in
$\cN$, respectively. We denote the predual of $\cN$ by $\cN_*$ and
its positive cone by $\cN_*^+$.
For a Hilbert space $\cK$ and $\xi\in\cK$, we define the vector functional
\[
\omega_\xi(x)\coloneqq\langle x\xi,\xi\rangle,
\qquad x\in\bB(\cK),
\]
and use the same symbol for its restriction to any von Neumann algebra acting on $\cK$.

We briefly recall the extended lower-semibounded self-adjoint part
of a von Neumann algebra. We refer to
\cite{haagerup_1979} and \cite{huw_2022} for details.

Following \cite[Definition~2.1]{huw_2022}, where the case of $\bB(\cH)$
is treated, we define the \emph{extended lower-semibounded self-adjoint
part} $\widehat{\cN}_{\mathrm{lb}}$ of an arbitrary von Neumann algebra
$\cN$ as the set of maps
\[
m\colon\cN_*^+\to(-\infty,+\infty]
\]
satisfying the following conditions:
\begin{enumerate}
\item
$m(\lambda\varphi)=\lambda m(\varphi)$
for $\varphi\in\cN_*^+$ and $\lambda\geq0$, where
$0\cdot(+\infty)=(+\infty)\cdot0=0$;
\item
$m(\varphi+\psi)=m(\varphi)+m(\psi)$
for $\varphi,\psi\in\cN_*^+$;
\item $m$ is lower semicontinuous on $\cN_*^+$ with respect to the
norm topology of $\cN_*$;
\item there exists $\ell\in\mathbb R$ such that
$m(\varphi)\geq\ell\varphi(1_{\cN})$
for $\varphi\in\cN_*^+$.
\end{enumerate}
The cone of all $m\in\widehat{\cN}_{\mathrm{lb}}$ satisfying
$m(\varphi)\geq0$ for every $\varphi\in\cN_*^+$ is Haagerup's
extended positive part $\widehat{\cN}_+$.

For $m,m_1,m_2\in\widehat{\cN}_{\mathrm{lb}}$,
$a\in\cN$, and $\lambda\geq0$, we define
\begin{align*}
(\lambda m)(\varphi)
&\coloneqq\lambda m(\varphi),\\
(m_1+m_2)(\varphi)
&\coloneqq m_1(\varphi)+m_2(\varphi),\\
(a^*ma)(\varphi)
&\coloneqq m(a\varphi a^*),
\end{align*}
where
$
(a\varphi a^*)(x)\coloneqq\varphi(a^*xa)
$ for $x\in\cN$.
We write $m_1\leq m_2$ if
$
m_1(\varphi)\leq m_2(\varphi)
$
for every $\varphi\in\cN_*^+$.
Every $a\in\cN_{\mathrm{sa}}$ is naturally identified with the
element of $\widehat{\cN}_{\mathrm{lb}}$ given by
$
\varphi\mapsto\varphi(a)$.
We say that $m\in\widehat{\cN}_{\mathrm{lb}}$ is \emph{bounded} if
$m$ belongs to the image of $\cN_{\mathrm{sa}}$ under this
identification, that is, if there exists $a\in\cN_{\mathrm{sa}}$
such that $m(\varphi)=\varphi(a)$ for every $\varphi\in\cN_*^+$.

We also use compression by operators between different Hilbert
spaces. If $\cH_1$ and $\cH_2$ are Hilbert spaces,
$C\in\bB(\cH_1,\cH_2)$, and
$m\in\widehat{\bB(\cH_2)}_{\mathrm{lb}}$, then
\[
(C^*mC)(\omega)\coloneqq m(C\omega C^*),
\qquad
\omega\in\bB(\cH_1)_*^+,
\]
where
\[
(C\omega C^*)(X)\coloneqq\omega(C^*XC),
\qquad X\in\bB(\cH_2).
\]
Then
\[
C^*mC\in\widehat{\bB(\cH_1)}_{\mathrm{lb}};
\]
see \cite[Lemma~2.4]{huw_2022}.

\begin{lemma}
\label{lem:vector-determined}
Let $m\in\widehat{\bB(\cH)}_{\mathrm{lb}}$, and let
$\omega\in\bB(\cH)_*^+$ be written as
\[
\omega=\sum_{n=1}^\infty\omega_{\xi_n},
\qquad
\sum_{n=1}^\infty\|\xi_n\|^2<+\infty.
\]
Then
\[
m(\omega)=\sum_{n=1}^\infty m(\omega_{\xi_n})
\in(-\infty,+\infty],
\]
where the series is well defined in $(-\infty,+\infty]$.
In particular, if $m_1,m_2\in\widehat{\bB(\cH)}_{\mathrm{lb}}$ satisfy
$m_1(\omega_\xi)=m_2(\omega_\xi)$ for every $\xi\in\cH$, then
$m_1=m_2$.
\end{lemma}

\begin{proof}
Let $\ell\in\mathbb R$ satisfy
$m(\varphi)\geq\ell\varphi(1)$ for $\varphi\in\bB(\cH)_*^+$.
Since
$m(\omega_{\xi_n})-\ell\|\xi_n\|^2\geq0$
and $\sum_n\|\xi_n\|^2<+\infty$, the series
$\sum_n m(\omega_{\xi_n})$ is well defined in $(-\infty,+\infty]$.
Set
\[
\omega_N\coloneqq\sum_{n=1}^N\omega_{\xi_n},
\qquad
r_N\coloneqq\omega-\omega_N\in\bB(\cH)_*^+.
\]
Then $\|r_N\|=r_N(1)=\sum_{n>N}\|\xi_n\|^2\to0$.
By additivity and the lower bound,
\[
m(\omega)
=
\sum_{n=1}^N m(\omega_{\xi_n})+m(r_N)
\geq
\sum_{n=1}^N m(\omega_{\xi_n})+\ell\,r_N(1),
\]
so that
$m(\omega)\geq\sum_{n=1}^\infty m(\omega_{\xi_n})$.
Conversely, since $\omega_N\to\omega$ in norm, the norm lower
semicontinuity of $m$ gives
\[
m(\omega)
\leq
\liminf_{N\to\infty}m(\omega_N)
=
\sum_{n=1}^\infty m(\omega_{\xi_n}).
\]
The last assertion follows because every
$\omega\in\bB(\cH)_*^+$ admits such a decomposition.
\end{proof}

The following spectral description is the von Neumann algebra
version of \cite[Proposition~2.2]{huw_2022}. 
It follows by applying \cite[Theorem~1.5]{haagerup_1979} to $m-\ell 1_{\cN}\in\widehat{\cN}_+$.

\begin{proposition}[cf.\ {\cite[Proposition 2.2]{huw_2022}}]
\label{prop:resolution}
Let $m\in\widehat{\cN}_{\mathrm{lb}}$, and let $\ell\in\mathbb R$
satisfy
\[
m(\varphi)\geq\ell\varphi(1_{\cN}),
\qquad \varphi\in\cN_*^+.
\]
Then there exist a unique projection $p_0\in\cN$ and a unique
nondecreasing right-continuous family of projections
$(e_t)_{t\in\mathbb R}$ in $\cN$ such that
$e_t=0$ for $t<\ell$,
\[
\lim_{t\to+\infty}e_t=p_0,
\]
and
\[
m(\varphi)
=
\int_{-\infty}^{+\infty}t\,d\varphi(e_t)
+
\infty\cdot\varphi(p_0^\perp),
\qquad
\varphi\in\cN_*^+.
\]
Moreover, $p_0^\perp=0$ if and only if
$
\{\varphi\in\cN_*^+\mid m(\varphi)<+\infty\}
$
is norm dense in $\cN_*^+$.
\end{proposition}

We call $p_0\cH$ and $p_0^\perp\cH$ the \emph{essential part} and
the \emph{infinite part} of $m$, respectively. We say that $m$ has
a \emph{dense domain} if $p_0^\perp=0$; see
\cite[Definition~2.3]{huw_2022}.

An element $m\in\widehat{\bB(\cH)}_{\mathrm{lb}}$ is said to be
\emph{affiliated with} $\cN$ if
$
u'^*mu'=m
$
for every $u'\in\cU(\cN')$.

As in Proposition~\ref{prop:resolution}, the following reduces to
\cite[Proposition~1.9]{haagerup_1979} by subtracting a scalar lower
bound $\ell 1$.

\begin{proposition}[cf.\ {\cite[Proposition~1.9]{haagerup_1979}}]
\label{prop:affiliated}
For $m\in\widehat{\cN}_{\mathrm{lb}}$, define
\[
\widetilde m(\omega)\coloneqq m(\omega|_{\cN}),
\qquad
\omega\in\bB(\cH)_*^+.
\]
Then the map $m\mapsto\widetilde m$ is a bijection between $\widehat{\cN}_{\mathrm{lb}}$ and the set of
elements of $\widehat{\bB(\cH)}_{\mathrm{lb}}$ affiliated with
$\cN$.
\end{proposition}

We also use the Pusz--Woronowicz (PW) functional calculus. 
We refer
to the original work \cite{pw}, the operator-theoretic reformulation
\cite{hatano_ueda_2021}, and the lower-semibounded extension
\cite{huw_2022}.

%%%%%%%%%%%%%%%%%%%%%%%%%%%%%%%%%%%%%%%%%%%%%%%%%%%%%%%%%%%%%%%%%%%%%%%%%%%%%%%%%%%%%%%%%%%%%%%%%%%%%%%%%%%%%%%%%%%%%%%%%%%%%%%%%%%%%%%%%%%%%%%%%%%%%%%%%%%%%%%%%%%%%%%%%%%%%%%%%%%%%%%%%%%%%%%%%%%%%%%%%%%%%%%%%%%%%%%%%%%%%%%%%%%%%%%%%%%%%%%%

\section{\texorpdfstring{Operator-valued maximal $f$-divergence}{Operator-valued maximal f-divergence}}

%%%%%%%%%%%%%%%%%%%%%%%%%%%%%%%%%%%%%%%%%%%%%%%%%%%%%%%%%%%%%%%%%%%%%%%%%%%%%%%%%%%%%%%%%%%%%%%%%%%%%%%%%%%%%%%%%%%%%%%%%%%%%%%%%%%%%%%%%%%%%%%%%%%%%%%%%%%%%%%%%%%%%%%%%%%%%%%%%%%%%%%%%%%%%%%%%%%%%%%%%%%%%%%%%%%%%%%%%%%%%%%%%%%%%%%%%%%%%%%%

Let $\cM$ and $\cN$ be von Neumann algebras, and assume that
$\cN$ acts on a Hilbert space $\cH$.
We denote by $\mathrm{CP}(\cM,\cN)$ the set of all completely
positive (CP) maps from $\cM$ to $\cN$.
Unless otherwise stated, maps in $\mathrm{CP}(\cM,\cN)$ are not
assumed to be normal.
For $\Phi, \Psi\in \mathrm{CP}(\cM,\cN)$,
we define $\Psi\leq_{\mathrm{cp}}\Phi$
if $\Phi-\Psi\in \mathrm{CP}(\cM,\cN)$.

We denote by $\mathrm{OC}(0,+\infty)$ the set of all real-valued
operator convex functions on $(0,+\infty)$.
Let $f\in\mathrm{OC}(0,+\infty)$.
Since $f$ is convex, there exist $a,b\in\mathbb R$ such that
\[
f(r)\geq ar+b,
\qquad r>0.
\]
Set
\[
f(+0)\coloneqq\lim_{t\downarrow0}f(t)\in(-\infty,+\infty],
\qquad
f'(+\infty)\coloneqq
\lim_{t\to+\infty}\frac{f(t)}{t}\in(-\infty,+\infty].
\]
Following \cite[(7.1), (7.2), and Definition~7.1]{huw_2022},
we define the perspective function
\[
\phi_f\colon\mathbb R_+^2\to(-\infty,+\infty]
\]
by
\[
\phi_f(s,t)\coloneqq
 \begin{cases}
  t f(s/t),&s,t>0,\\
  f'(+\infty)s,&s>0,\ t=0,\\
  f(+0)t,&s=0,\ t>0,\\
  0,&s=t=0.
 \end{cases}
\]
with the convention $(+\infty)\cdot 0=0$.
Then $\phi_f$ is a Borel function that is homogeneous and locally
bounded from below on $\mathbb R_+^2$. Hence, by the PW functional calculus
\cite[Definition~4.1, Theorem~4.3, and Definition~7.1]{huw_2022},
\[
\phi_f(A,B)\in\widehat{\bB(\cK)}_{\mathrm{lb}}
\]
is defined for every Hilbert space $\cK$ and all
$A,B\in\bB(\cK)_+$.

We also use the PW functional calculus for positive
sesquilinear forms.

\begin{definition}
\label{def:extended-pw-form}
Let $\cV$ be a complex vector space, and let $\alpha$ and $\beta$
be positive sesquilinear forms on $\cV$. 
Let $\cH_\gamma$ be the Hilbert space obtained by separation and
completion with respect to $\gamma\coloneqq\alpha+\beta$, and let
$
h_\gamma\colon\cV\to\cH_\gamma
$
be the canonical linear map.
Then there exist unique positive contractions
\[
R_{\alpha,\beta},S_{\alpha,\beta}
\in\bB(\cH_\gamma)_+
\]
such that
\[
\alpha(X,Y)
=
\langle
R_{\alpha,\beta}h_\gamma(X),h_\gamma(Y)
\rangle,
\qquad
\beta(X,Y)
=
\langle
S_{\alpha,\beta}h_\gamma(X),h_\gamma(Y)
\rangle,
\]
and
\[
R_{\alpha,\beta}+S_{\alpha,\beta}
=
I_{\cH_\gamma}.
\]

We define the lower semibounded quadratic form
\[
q_f[\alpha,\beta]\colon\cV\to(-\infty,+\infty]
\]
by
\[
q_f[\alpha,\beta](X)
\coloneqq
\phi_f(R_{\alpha,\beta},S_{\alpha,\beta})
\bigl(\omega_{h_\gamma(X)}\bigr),
\qquad X\in\cV.
\]
\end{definition}

\blem
\label{lem:compatible-form}
Let $\alpha$ and $\beta$ be positive sesquilinear forms on $\cV$,
and let
\[
(h,A,B,\cK)
\]
be a compatible representation of $(\alpha,\beta)$,
i.e.,
$h\colon\cV\to\cK$ has dense range,
\[
A,B\in\bB(\cK)_+,\qquad
AB=BA,
\]
and
\[
\alpha(X,Y)=\langle Ah(X),h(Y)\rangle,
\qquad
\beta(X,Y)=\langle Bh(X),h(Y)\rangle.
\]
Then
\[
q_f[\alpha,\beta](X)
=
\phi_f(A,B)(\omega_{h(X)}),
\qquad X\in\cV.
\]
\elem

\bpf
Both $(h_\gamma,R_{\alpha,\beta},S_{\alpha,\beta},\cH_\gamma)$ and
$(h,A,B,\cK)$ are compatible representations of $(\alpha,\beta)$,
and in each of them the two operators commute. 
We simply write $R\coloneqq R_{\alpha,\beta}$
and $S\coloneqq S_{\alpha,\beta}$.

For a commuting pair
$(C,D)$ of positive operators we write $E_{C,D}$ for its joint
spectral measure on $\mathbb R_+^2$, so that
\[
\phi_f(C,D)(\omega_\zeta)
=
\int_{\mathbb R_+^2}\phi_f(s,t)\,d\langle E_{C,D}(s,t)\zeta,\zeta\rangle
\]
by \cite[Definition~4.1(1) and Theorem~4.3]{huw_2022}.

Choose $a,b\in\mathbb R$ with $f(r)\geq ar+b$ for $r>0$. Then
\[
\phi_f(s,t)\geq as+bt,
\qquad (s,t)\in\mathbb R_+^2,
\]
and we decompose
\[
\phi_f(s,t)
=
as+bt
+
g(s,t)\,\mathbf 1_{\{s,t>0\}}
+
\bigl(f'(+\infty)-a\bigr)s\,\mathbf 1_{\{t=0\}}
+
\bigl(f(+0)-b\bigr)t\,\mathbf 1_{\{s=0\}},
\]
where $g(s,t)\coloneqq\phi_f(s,t)-as-bt\geq0$.

For $n\in\mathbb N$, let
\[
\mathcal C_n\coloneqq\{(s,t)\in\mathbb R_+^2: s/n\leq t\leq ns\}.
\]
Since $g$ is continuous on the compact set
$\mathcal C_n\cap\{s+t=1\}$,
$g_n\coloneqq g\,\mathbf 1_{\mathcal C_n}$ is a nonnegative
homogeneous Borel function satisfying
$g_n(s,t)\leq M_n(s+t)$ for some $M_n<\infty$.
Moreover, $g_n\uparrow g\,\mathbf 1_{\{s,t>0\}}$ pointwise.

The functions $as+bt$, $g_n$, $s\,\mathbf 1_{\{t=0\}}$ and
$t\,\mathbf 1_{\{s=0\}}$ are homogeneous Borel functions that are
bounded on $\{s+t=1\}$. Hence, by \cite[Theorem~1.2]{pw},
\begin{equation}\label{eq:pw-bounded}
u(R,S)(\omega_{h_\gamma(X)})
=
u(A,B)(\omega_{h(X)})
\end{equation}
for
\[
u\in\bigl\{\,as+bt,\ g_n,\ s\mathbf 1_{\{t=0\}},\ t\mathbf 1_{\{s=0\}}\,\bigr\}.
\]
Since $0\leq g_n\uparrow g\mathbf 1_{\{s,t>0\}}$, the monotone
convergence theorem and \eqref{eq:pw-bounded} give
\begin{align*}
(g\mathbf 1_{\{s,t>0\}})(R,S)(\omega_{h_\gamma(X)})
&=
\lim_{n\to\infty}g_n(R,S)(\omega_{h_\gamma(X)})
\\
&=
\lim_{n\to\infty}g_n(A,B)(\omega_{h(X)})
\\
&=
(g\mathbf 1_{\{s,t>0\}})(A,B)(\omega_{h(X)})
\in[0,+\infty].
\end{align*}
Moreover, by \eqref{eq:pw-bounded} and the convention
$(+\infty)\cdot0=0$,
\begin{align*}
&\bigl(f'(+\infty)-a\bigr)\,
(s\mathbf 1_{\{t=0\}})(R,S)(\omega_{h_\gamma(X)})
\\
&\qquad\qquad=
\bigl(f'(+\infty)-a\bigr)\,
(s\mathbf 1_{\{t=0\}})(A,B)(\omega_{h(X)})
\in[0,+\infty],\\
&\bigl(f(+0)-b\bigr)\,
(t\mathbf 1_{\{s=0\}})(R,S)(\omega_{h_\gamma(X)})
\\
&\qquad\qquad=
\bigl(f(+0)-b\bigr)\,
(t\mathbf 1_{\{s=0\}})(A,B)(\omega_{h(X)})
\in[0,+\infty].
\end{align*}
Summing the four terms gives
\[
\phi_f(R,S)(\omega_{h_\gamma(X)})
=
\phi_f(A,B)(\omega_{h(X)}),
\]
that is,
$q_f[\alpha,\beta](X)=\phi_f(A,B)(\omega_{h(X)})$.
\epf

For arbitrary $A,B\in\bB(\cH)_+$, let
\[
(T_{A,B},R_{A,B},S_{A,B},\cH_{A,B})
\]
be their canonical compatible representation given in
\cite[Lemma~4.4]{huw_2022}.
Set
\[
s_A(\xi,\eta)
=
\langle
A\xi,\eta
\rangle,
\qquad
s_B(\xi,\eta)
=
\langle
B\xi,\eta
\rangle.
\]
Applying
Lemma~\ref{lem:compatible-form} to this representation,
we obtain
\[
q_f[s_A,s_B](\xi)
=
\phi_f(R_{A,B},S_{A,B})
\bigl(\omega_{T_{A,B}\xi}\bigr)
=
\phi_f(A,B)(\omega_\xi),
\qquad \xi\in\cH.
\]
See \cite[Lemma~4.4, Remark~4.5, and (4.4)]{huw_2022}.

%%%%%%%%%%%%%%%%%%%%%%%%%%%%%%%%%%%%%%%%%%%%%%%%%%%%%%%%%%%%%%%%%%%%%%%%%%%%%%%%%%%%%%%%%%%%%%%%%%%%%%%%%%%%%%%%%%%%%%%%%%%%%%%%%%%%%%%%%%%%%%%%%%%%%%%%%%%%%%%%%%%%%%%%%%%%%%%%%%%%%%%%%%%%%%%%%%%%%%%%%%%%%%%%%%%%%%%%%%%%%%%%%%%%%%%%%%%%%%%%

%%%%%%%%%%%%%%%%%%%%%%%%%%%%%%%%%%%%%%%%%%%%%%%%%%%%%%%%%%%%%%%%%%%%%%%%%%%%%%%%%%%%%%%%%%%%%%%%%%%%%%%%%%%%%%%%%%%%%%%%%%%%%%%%%%%%%%%%%%%%%%%%%%%%%%%%%%%%%%%%%%%%%%%%%%%%%%%%%%%%%%%%%%%%%%%%%%%%%%%%%%%%%%%%%%%%%%%%%%%%%%%%%%%%%%%%%%%%%%%%

\bdf
To simplify notation, we adopt the following convention.
Let $\Phi,\Psi,\Gamma\in\mathrm{CP}(\cM,\cN)$ satisfy
\[
\Phi,\Psi\leq_{\mathrm{cp}}\Gamma.
\]
A \emph{minimal Stinespring--Radon--Nikodym (SRN) realization of
$(\Phi,\Psi)$ relative to} $\Gamma$ is a quintuple
\[
(\pi,V,\cK;A',B')
\]
such that $(\pi,V,\cK)$ is a minimal Stinespring representation
of $\Gamma$ and $A',B'\in\pi(\cM)'$ are the corresponding RN
derivatives. Thus,
\[
0\leq A',B'\leq I_{\cK}
\]
and
\[
\Phi(x)=V^*A'\pi(x)V,
\qquad
\Psi(x)=V^*B'\pi(x)V,
\qquad x\in\cM.
\]
The existence of such a realization follows from Stinespring's
theorem and the Radon--Nikodym theorem for CP maps; see
\cite{stinespring_1955,arveson_1969,belavkin-staszewski_1986}.

Unless otherwise specified, 
we take $\Gamma=\Phi+\Psi$ and simply
call $(\pi,V,\cK;A',B')$ a \emph{minimal SRN realization} of
$(\Phi,\Psi)$. 
In this case, uniqueness of the RN derivative of
$\Gamma$ gives
\[
A'+B'=I_{\cK}.
\]
\edf

\begin{definition}
Let $f\in\mathrm{OC}(0,+\infty)$, and let
$(\pi,V,\cK;A',B')$ be a minimal SRN realization of
$(\Phi,\Psi)$. The \emph{maximal $f$-divergence} of $\Phi$
relative to $\Psi$ is defined by
\[
\widehat S_f(\Phi\|\Psi)
\coloneqq
V^*\phi_f(A',B')V
\in
\widehat{\bB(\cH)}_{\mathrm{lb}}.
\]
\end{definition}

The right-hand side is independent of the choice of a minimal SRN realization. 
Indeed, this follows from the unitary uniqueness of
minimal Stinespring representations, 
the uniqueness of the RN derivatives, 
and the unitary covariance of the PW functional calculus. More generally, independence of the dominating CP map
is established in Proposition~\ref{prop:independent}.

We first prove that $\widehat S_f(\Phi\|\Psi)$
indeed belongs to $\widehat{\cN}_{\mathrm{lb}}$.

\blem\label{lem:pi-commutant}
There exists a unital $*$-representation
$\pi'\colon\cN'\to\bB(\cK)$ satisfying
\[
\pi'(b')
\left(\sum_i\pi(x_i)V\xi_i\right)
=
\sum_i\pi(x_i)V(b'\xi_i)
\]
for $x_i\in\cM$, $\xi_i\in\cH$, and $b'\in\cN'$.
In particular, $\pi'(b')V=Vb'$.

Moreover,
\[
\pi'(\cN')\subseteq\pi(\cM)'
\]
and
\[
A', B'\in\pi'(\cN')'.
\]
\elem

\bpf
Let $b'\in\cN'$ and
\[
\zeta=\sum_{i=1}^n\pi(x_i)V\xi_i.
\]
Set
\[
\boldsymbol{\xi}=(\xi_1,\ldots,\xi_n)^T,\qquad
G=[\Gamma(x_i^*x_j)]_{i,j}\in M_n(\cN)_+,
\]
and $D=\operatorname{diag}(b',\ldots,b')$.
Since $G$ and $D$ commute, we have
\begin{align*}
\left\|
\sum_i\pi(x_i)V(b'\xi_i)
\right\|^2
&=
\langle D^*GD\boldsymbol{\xi},\boldsymbol{\xi}\rangle\\
&=
\langle
G^{1/2}D^*DG^{1/2}\boldsymbol{\xi},
\boldsymbol{\xi}
\rangle\\
&\leq
\|b'\|^2
\langle G\boldsymbol{\xi},\boldsymbol{\xi}\rangle
=
\|b'\|^2\|\zeta\|^2.
\end{align*}
Thus the stated formula defines a bounded operator $\pi'(b')$
on $\cK$.
Moreover $\pi'$ is a
unital $*$-representation. The definition also gives
$
\pi'(b')\pi(a)=\pi(a)\pi'(b')
$ for $a\in\cM$,
so that $\pi'(\cN')\subseteq\pi(\cM)'$.

Finally, let $u'\in\cU(\cN')$.
Since $\Phi(x)\in\cN$
and
$\pi'(u')\pi(x)=\pi(x)\pi'(u')$,
\[
\Phi(x)
=
u'^*\Phi(x)u'
=
u'^*V^*A'\pi(x)Vu'
=
V^*\pi'(u')^*A'\pi'(u')\pi(x)V.
\]
The uniqueness of RN derivatives implies
$
\pi'(u')^*A'\pi'(u')=A'.
$
Thus $A'$ commutes with $\pi'(u')$ for every
$u'\in\cU(\cN')$, and hence $A'\in\pi'(\cN')'$.
Similarly, $B'\in\pi'(\cN')'$.
\epf

\blem\label{lem:N-valued}
The maximal $f$-divergence $\widehat S_f(\Phi\|\Psi)$ is
affiliated with $\cN$. Consequently,
\[
\widehat S_f(\Phi\|\Psi)
\in\widehat{\cN}_{\mathrm{lb}}.
\]
\elem

\bpf
Let $u'\in\cU(\cN')$ and put $U'\coloneqq\pi'(u')$.
By Lemma~\ref{lem:pi-commutant}, $U'$ commutes with $A'$ and $B'$.
Hence the unitary covariance of the PW functional calculus gives
\[
U'^*\phi_f(A',B')U'=\phi_f(A',B').
\]
Since $U'V=Vu'$, we obtain
\[
u'^*\widehat S_f(\Phi\|\Psi)u'
=
V^*U'^*\phi_f(A',B')U'V
=
\widehat S_f(\Phi\|\Psi).
\]
Thus $\widehat S_f(\Phi\|\Psi)$ is affiliated with $\cN$, and
the conclusion follows from Proposition~\ref{prop:affiliated}.
\epf

We next show that the maximal $f$-divergence 
$\widehat S_f(\Phi\|\Psi)$ can be computed using
any CP map that dominates $\Phi$ and $\Psi$ up to positive constants.

\bprop
\label{prop:independent}
Let $\Gamma\in\mathrm{CP}(\cM,\cN)$ satisfy
\[
\Phi\leq_{\mathrm{cp}}c_\Phi\Gamma,
\qquad
\Psi\leq_{\mathrm{cp}}c_\Psi\Gamma
\]
for some $c_\Phi,c_\Psi>0$. Let $(\pi,V,\cK)$ be a minimal
Stinespring representation of $\Gamma$, and let
$A',B'\in\pi(\cM)'_+$ be the corresponding RN derivatives:
\[
\Phi(x)=V^*A'\pi(x)V,
\qquad
\Psi(x)=V^*B'\pi(x)V,
\qquad x\in\cM.
\]
Then
\[
\widehat S_f(\Phi\|\Psi)
=
V^*\phi_f(A',B')V.
\]
In particular, the right-hand side is independent of the choice
of the dominating CP map $\Gamma$.
\eprop

\bpf
For $\Theta\in\mathrm{CP}(\cM,\cN)$, define the positive
sesquilinear form $s_\Theta$ on $\cM\odot\cH$ by
\[
s_\Theta
\left(
\sum_i x_i\otimes\xi_i,
\sum_j y_j\otimes\eta_j
\right)
\coloneqq
\sum_{i,j}
\left\langle
\Theta(y_j^*x_i)\xi_i,\eta_j
\right\rangle.
\]
Let
$
k_\Gamma\colon\cM\odot\cH\to\cK
$
be the canonical map given by
\[
k_\Gamma
\left(
\sum_i x_i\otimes\xi_i
\right)
\coloneqq
\sum_i\pi(x_i)V\xi_i.
\]
The minimality of $(\pi,V,\cK)$ implies that
$k_\Gamma(\cM\odot\cH)$ is dense in $\cK$. Moreover,
\[
s_\Phi(X,Y)
=
\langle A'k_\Gamma(X),k_\Gamma(Y)\rangle,
\qquad
s_\Psi(X,Y)
=
\langle B'k_\Gamma(X),k_\Gamma(Y)\rangle.
\]

Let
\[
(T_{A',B'},R_{A',B'},S_{A',B'},\cH_{A',B'})
\]
be the canonical compatible representation of $(A',B')$
given in \cite[Lemma~4.4]{huw_2022}. Thus,
\[
A'
=
T_{A',B'}^*R_{A',B'}T_{A',B'},
\qquad
B'
=
T_{A',B'}^*S_{A',B'}T_{A',B'},
\]
and
\[
R_{A',B'}+S_{A',B'}
=
I_{\cH_{A',B'}}.
\]
Set
\[
h_\Gamma
\coloneqq
T_{A',B'}k_\Gamma
\colon
\cM\odot\cH\to\cH_{A',B'}.
\]
Since $k_\Gamma(\cM\odot\cH)$ is dense in $\cK$ and
$T_{A',B'}$ has dense range, we have
\[
\overline{h_\Gamma(\cM\odot\cH)}
=
\cH_{A',B'}.
\]
Moreover,
\[
s_\Phi(X,Y)
=
\langle R_{A',B'}h_\Gamma(X),h_\Gamma(Y)\rangle,
\qquad
s_\Psi(X,Y)
=
\langle S_{A',B'}h_\Gamma(X),h_\Gamma(Y)\rangle.
\]
Thus
\[
(h_\Gamma,R_{A',B'},S_{A',B'},\cH_{A',B'})
\]
is a compatible representation of $(s_\Phi,s_\Psi)$. Hence,
Lemma~\ref{lem:compatible-form} and the canonical realization
\cite[(4.4)]{huw_2022} give
\begin{align}
q_f[s_\Phi,s_\Psi](X)
&=
\phi_f(R_{A',B'},S_{A',B'})
\bigl(\omega_{h_\Gamma(X)}\bigr)
\nonumber\\
&=
\phi_f(R_{A',B'},S_{A',B'})
\bigl(
\omega_{T_{A',B'}k_\Gamma(X)}
\bigr)
\nonumber\\
&=
\bigl(
T_{A',B'}^*
\phi_f(R_{A',B'},S_{A',B'})
T_{A',B'}
\bigr)
\bigl(\omega_{k_\Gamma(X)}\bigr)
\nonumber\\
&=
\phi_f(A',B')
\bigl(\omega_{k_\Gamma(X)}\bigr).
\label{eq:form-realization}
\end{align}

Let
\[
(\pi_0,V_0,\cK_0;A_0',B_0')
\]
be the minimal SRN realization associated with
$\Gamma_0=\Phi+\Psi$ that is used in the definition of
$\widehat S_f(\Phi\|\Psi)$. 
By the same argument applied to $\Gamma_0$, we have
\[
q_f[s_\Phi,s_\Psi](X)
=
\phi_f(A_0',B_0')
\bigl(\omega_{k_0(X)}\bigr),
\]
where
\[
k_0
\left(
\sum_i x_i\otimes\xi_i
\right)
=
\sum_i\pi_0(x_i)V_0\xi_i.
\]

For an arbitrary $\xi\in\cH$, take
\[
X=1_{\cM}\otimes\xi.
\]
Then we have
\[
k_\Gamma(X)=V\xi,
\qquad
k_0(X)=V_0\xi.
\]
Consequently,
\begin{align*}
\bigl(V^*\phi_f(A',B')V\bigr)(\omega_\xi)
&=
q_f[s_\Phi,s_\Psi](1_{\cM}\otimes\xi)
\\
&=
\bigl(V_0^*\phi_f(A_0',B_0')V_0\bigr)(\omega_\xi)
=
\widehat S_f(\Phi\|\Psi)(\omega_\xi).
\end{align*}

Since $\xi\in\cH$ is arbitrary, Lemma~\ref{lem:vector-determined},
applied to the two elements
$V^*\phi_f(A',B')V$ and $V_0^*\phi_f(A_0',B_0')V_0$ of
$\widehat{\bB(\cH)}_{\mathrm{lb}}$, yields
\[
V^*\phi_f(A',B')V
=
V_0^*\phi_f(A_0',B_0')V_0
=
\widehat S_f(\Phi\|\Psi)
\]
in $\widehat{\bB(\cH)}_{\mathrm{lb}}$, and hence also in
$\widehat{\cN}_{\mathrm{lb}}$ by Lemma~\ref{lem:N-valued} and
Proposition~\ref{prop:affiliated}.
\epf

We next compute the maximal $f$-divergence 
$\widehat S_f(\Phi\|\Psi)$ using a Stinespring
representation of $\Phi+\Psi$ that need not be minimal.

\blem\label{lem:isometry}
Let $\cH$ and $\cK$ be Hilbert spaces, let
$C,D\in\bB(\cK)_+$, and let $U\colon\cH\to\cK$ be an isometry.
If $C$ and $D$ commute with $P\coloneqq UU^*$, then
\[
\phi_f(U^*CU,U^*DU)
=
U^*\phi_f(C,D)U.
\]
\elem

\bpf
With respect to the decomposition
\[
\cK=P\cK\oplus P^\perp\cK,
\]
write
\[
C=C_0\oplus C_1,
\qquad
D=D_0\oplus D_1.
\]
Regarded as a map from $\cH$ onto $P\cK$, the operator $U$ is
unitary. Hence, by unitary covariance and the direct sum property of
the PW functional calculus
\cite[Proposition~4.6]{huw_2022},
\[
\phi_f(U^*CU,U^*DU)
=
\phi_f(U^*C_0U,U^*D_0U)
=
U^*\phi_f(C_0,D_0)U
=
U^*\phi_f(C,D)U.
\]
\epf

Let $\Phi,\Psi\in\mathrm{CP}(\cM,\cN)$ and set
$
\Gamma\coloneqq\Phi+\Psi.
$
Let $(\pi_1,V_1,\cK_1)$ be a Stinespring representation of
$\Gamma$ that is not assumed to be minimal, and let
$A_1',B_1'\in\pi_1(\cM)'_+$ be RN derivatives of $\Phi$ and
$\Psi$, respectively.
Thus,
\[
\Phi(a)=V_1^*A_1'\pi_1(a)V_1,
\qquad
\Psi(a)=V_1^*B_1'\pi_1(a)V_1,
\qquad a\in\cM.
\]
The RN derivatives need not be unique because the Stinespring
representation is not assumed to be minimal.

Let
$
(\pi_0,V_0,\cK_0;A_0',B_0')
$
be a minimal SRN realization of $(\Phi,\Psi)$. By the standard
uniqueness theorem for Stinespring representations, there exists a
unique isometry
$
W\colon\cK_0\to\cK_1
$
such that
\[
WV_0=V_1,
\qquad
W\pi_0(a)=\pi_1(a)W,
\qquad a\in\cM.
\]
Explicitly, $W$ is determined on the canonical dense subspace by
\[
W\left(\sum_i\pi_0(x_i)V_0\xi_i\right)
=
\sum_i\pi_1(x_i)V_1\xi_i.
\]
Set
$
P\coloneqq WW^*.
$
Then
\[
P\cK_1
=
\overline{\pi_1(\cM)V_1\cH},
\qquad
P\in\pi_1(\cM)'.
\]

\bprop\label{prop:non-minimal}
For every choice of RN derivatives $A_1'$ and $B_1'$ as above,
\[
\widehat S_f(\Phi\|\Psi)
\leq
V_1^*\phi_f(A_1',B_1')V_1.
\]
If $A_1'$ and $B_1'$ commute with $P=WW^*$, then equality holds:
\[
\widehat S_f(\Phi\|\Psi)
=
V_1^*\phi_f(A_1',B_1')V_1.
\]
\eprop

\bpf
The intertwining relations imply that
\[
W^*A_1'W,\ W^*B_1'W\in\pi_0(\cM)',
\]
and
\[
V_0^*(W^*A_1'W)\pi_0(a)V_0=\Phi(a),
\qquad
V_0^*(W^*B_1'W)\pi_0(a)V_0=\Psi(a).
\]
By uniqueness of the RN derivatives in the minimal Stinespring
representation,
\[
W^*A_1'W=A_0',
\qquad
W^*B_1'W=B_0'.
\]

By the isometric compression inequality for the operator perspective
\cite[Theorem~4.10(iii) and Theorem~7.2(1)]{huw_2022},
we obtain
\[
\phi_f(A_0',B_0')
=
\phi_f(W^*A_1'W,W^*B_1'W)
\leq
W^*\phi_f(A_1',B_1')W.
\]
Consequently,
\[
\widehat S_f(\Phi\|\Psi)
=
V_0^*\phi_f(A_0',B_0')V_0
\leq
V_0^*W^*\phi_f(A_1',B_1')WV_0
=
V_1^*\phi_f(A_1',B_1')V_1.
\]

If $A_1'$ and $B_1'$ commute with $P$, then Lemma~\ref{lem:isometry}
yields
\[
\phi_f(W^*A_1'W,W^*B_1'W)
=
W^*\phi_f(A_1',B_1')W,
\]
and hence equality follows.
\epf

%%%%%%%%%%%%%%%%%%%%%%%%%%%%%%%%%%%%%%%%%%%%%%%%%%%%%%%%%%%%%%%%%%%%%%%%%%%%%%%%%%%%%%%%%%%%%%%%%%%%%%%%%%%%%%%%%%%%%%%%%%%%%%%%%%%%%%%%%%%%%%%%%%%%%%%%%%%%%%%%%%%%%%%%%%%%%%%%%%%%%%%%%%%%%%%%%%%%%%%%%%%%%%%%%%%%%%%%%%%%%%%%%%%%%%%%%%%%%%%%%%%%%%%%%%%%%%%%%%%%%
\section{Fundamental properties}
%%%%%%%%%%%%%%%%%%%%%%%%%%%%%%%%%%%%%%%%%%%%%%%%%%%%%%%%%%%%%%%%%%%%%%%%%%%%%%%%%%%%%%%%%%%%%%%%%%%%%%%%%%%%%%%%%%%%%%%%%%%%%%%%%%%%%%%%%%%%%%%%%%%%%%%%%%%%%%%%%%%%%%%%%%%%%%%%%%%%%%%%%%%%%%%%%%%%%%%%%%%%%%%%%%%%%%%%%%%%%%%%%%%%%%%%%%%%%%%%%%%%%%%%%%%%%%%%%%%%%

Throughout this section we fix
$f\in\mathrm{OC}(0,+\infty)$.
We establish several fundamental properties of the maximal
$f$-divergence for CP maps.

%%%%%%%%%%%%%%%%%%%%%%%%%%%%%%%%%%%%%%%%%%%%%%%%%%%%%%%%%%%%%%%%%%%%%%%%%%%%%%%%%%%%%%%%%%%%%%%%%%%%%%%%%%%%%%%%%%%%%%%%%%%%%%%%%%%%%%%%%%%%%%%%%%%%%%%%%%%%%%%%%%%%%%%%%%%%%%%%%%%%%%%%%%%%%%%%%%%%%%%%%%%%%%%%%%%%%%%%%%%%%%%%%%%%%%%%%%%%%%%%%%%%%%%%%%%%%%%%%%%%%
\subsection{Joint subadditivity}
%%%%%%%%%%%%%%%%%%%%%%%%%%%%%%%%%%%%%%%%%%%%%%%%%%%%%%%%%%%%%%%%%%%%%%%%%%%%%%%%%%%%%%%%%%%%%%%%%%%%%%%%%%%%%%%%%%%%%%%%%%%%%%%%%%%%%%%%%%%%%%%%%%%%%%%%%%%%%%%%%%%%%%%%%%%%%%%%%%%%%%%%%%%%%%%%%%%%%%%%%%%%%%%%%%%%%%%%%%%%%%%%%%%%%%%%%%%%%%%%%%%%%%%%%%%%%%%%%%%%%

\bthm\label{thm:subadditivity}
Let $n\in\mathbb N$ and let
$\Phi_i,\Psi_i\in\mathrm{CP}(\cM,\cN)$ for $1\leq i\leq n$.
Then
\[
\widehat S_f
\left(
\sum_{i=1}^n\Phi_i
\,\middle\|\,
\sum_{i=1}^n\Psi_i
\right)
\leq
\sum_{i=1}^n\widehat S_f(\Phi_i\|\Psi_i).
\]
\ethm

\bpf
Set
\[
\Gamma\coloneqq\sum_{i=1}^n(\Phi_i+\Psi_i),
\]
and let $(\pi,V,\cK)$ be a minimal Stinespring representation of
$\Gamma$. For each $i$, let $A_i',B_i'\in\pi(\cM)'_+$ be the RN
derivatives of $\Phi_i,\Psi_i$ with respect to $\Gamma$.
Then
\[
\sum_{i=1}^nA_i'
\quad\text{and}\quad
\sum_{i=1}^nB_i'
\]
are the RN derivatives of $\sum_i\Phi_i$ and $\sum_i\Psi_i$,
respectively, and
\[
\sum_{i=1}^n(A_i'+B_i')=I_{\cK}.
\]
By Proposition~\ref{prop:independent} and the joint subadditivity of the
operator perspective
\cite[Theorem~7.2(1)]{huw_2022}, we obtain
\begin{align*}
\widehat S_f
\left(
\sum_{i=1}^n\Phi_i
\,\middle\|\,
\sum_{i=1}^n\Psi_i
\right)
&=
V^*\phi_f
\left(
\sum_{i=1}^nA_i',
\sum_{i=1}^nB_i'
\right)V\\
&\leq
\sum_{i=1}^nV^*\phi_f(A_i',B_i')V
=
\sum_{i=1}^n\widehat S_f(\Phi_i\|\Psi_i).
\end{align*}
\epf

%%%%%%%%%%%%%%%%%%%%%%%%%%%%%%%%%%%%%%%%%%%%%%%%%%%%%%%%%%%%%%%%%%%%%%%%%%%%%%%%%%%%%%%%%%%%%%%%%%%%%%%%%%%%%%%%%%%%%%%%%%%%%%%%%%%%%%%%%%%%%%%%%%%%%%%%%%%%%%%%%%%%%%%%%%%%%%%%%%%%%%%%%%%%%%%%%%%%%%%%%%%%%%%%%%%%%%%%%%%%%%%%%%%%%%%%%%%%%%%%%%%%%%%%%%%%%%%%%%%%%
\subsection{Transpose}
%%%%%%%%%%%%%%%%%%%%%%%%%%%%%%%%%%%%%%%%%%%%%%%%%%%%%%%%%%%%%%%%%%%%%%%%%%%%%%%%%%%%%%%%%%%%%%%%%%%%%%%%%%%%%%%%%%%%%%%%%%%%%%%%%%%%%%%%%%%%%%%%%%%%%%%%%%%%%%%%%%%%%%%%%%%%%%%%%%%%%%%%%%%%%%%%%%%%%%%%%%%%%%%%%%%%%%%%%%%%%%%%%%%%%%%%%%%%%%%%%%%%%%%%%%%%%%%%%%%%%

The \emph{transpose} of $f$ is defined by
\[
\widetilde f(t)\coloneqq t f(t^{-1}),
\qquad t>0.
\]
Then $\widetilde f\in\mathrm{OC}(0,+\infty)$ and
\[
\widetilde f'(+\infty)=f(+0),
\qquad
\widetilde f(+0)=f'(+\infty).
\]
Consequently,
\[
\phi_{\widetilde f}(s,t)=\phi_f(t,s),
\qquad (s,t)\in\mathbb R_+^2.
\]
See \cite[Definition~7.1, Theorem~7.2(1), and (7.3)]{huw_2022}.

\bprop\label{prop:transpose}
Let $\Phi,\Psi\in\mathrm{CP}(\cM,\cN)$. Then
\[
\widehat S_{\widetilde f}(\Phi\|\Psi)
=
\widehat S_f(\Psi\|\Phi).
\]
\eprop

\bpf
Let $(\pi,V,\cK;A',B')$ be a minimal SRN realization of
$(\Phi,\Psi)$. Since $(\pi,V,\cK;B',A')$ is a minimal SRN
realization of $(\Psi,\Phi)$, we have
\[
\widehat S_{\widetilde f}(\Phi\|\Psi)
=
V^*\phi_{\widetilde f}(A',B')V
=
V^*\phi_f(B',A')V
=
\widehat S_f(\Psi\|\Phi).
\]
\epf

%%%%%%%%%%%%%%%%%%%%%%%%%%%%%%%%%%%%%%%%%%%%%%%%%%%%%%%%%%%%%%%%%%%%%%%%%%%%%%%%%%%%%%%%%%%%%%%%%%%%%%%%%%%%%%%%%%%%%%%%%%%%%%%%%%%%%%%%%%%%%%%%%%%%%%%%%%%%%%%%%%%%%%%%%%%%%%%%%%%%%%%%%%%%%%%%%%%%%%%%%%%%%%%%%%%%%%%%%%%%%%%%%%%%%%%%%%%%%%%%%%%%%%%%%%%%%%%%%%%%%
\subsection{Precomposition monotonicity}
%%%%%%%%%%%%%%%%%%%%%%%%%%%%%%%%%%%%%%%%%%%%%%%%%%%%%%%%%%%%%%%%%%%%%%%%%%%%%%%%%%%%%%%%%%%%%%%%%%%%%%%%%%%%%%%%%%%%%%%%%%%%%%%%%%%%%%%%%%%%%%%%%%%%%%%%%%%%%%%%%%%%%%%%%%%%%%%%%%%%%%%%%%%%%%%%%%%%%%%%%%%%%%%%%%%%%%%%%%%%%%%%%%%%%%%%%%%%%%%%%%%%%%%%%%%%%%%%%%%%%

\bthm\label{thm:precomposition}
Let $\Phi,\Psi\in\mathrm{CP}(\cM,\cN)$, let $\cL$ be a von
Neumann algebra, and let
$\Theta\in\mathrm{CP}(\cL,\cM)$ be unital. Then
\[
\widehat S_f(\Phi\circ\Theta\|\Psi\circ\Theta)
\leq
\widehat S_f(\Phi\|\Psi).
\]
\ethm

\bpf
Let
$(\pi,V,\cK;A',B')$
be a minimal SRN realization of $(\Phi,\Psi)$, and set
$\Gamma\coloneqq\Phi+\Psi$.
Thus
\[
A'+B'=I_{\cK}.
\]
Since
$
\pi\circ\Theta\colon\cL\to\bB(\cK)
$
is a unital CP map, we take
a minimal Stinespring
representation $(\widetilde\pi,W,\widetilde\cK)$ of $\pi\circ\Theta$
such that
$
W\colon\cK\to\widetilde\cK
$
is an isometry.

Set
\[
\widetilde V\coloneqq WV\colon\cH\to\widetilde\cK.
\]
Then, for every $a\in\cL$,
\[
(\Gamma\circ\Theta)(a)
=
V^*\pi(\Theta(a))V
=
V^*W^*\widetilde\pi(a)WV
=
\widetilde V^*\widetilde\pi(a)\widetilde V.
\]
Thus $(\widetilde\pi,\widetilde V,\widetilde\cK)$
is a Stinespring representation of $\Gamma\circ\Theta$ that need not be minimal.

We next construct the required RN derivatives. Set
\[
\mathcal C\coloneqq(\pi\circ\Theta)(\cL)'
\subseteq\bB(\cK).
\]
By the same arguments as in the proof of Lemma~\ref{lem:pi-commutant},
we obtain a unital $*$-representation
\[
\rho\colon\mathcal C\to\bB(\widetilde\cK),
\]
defined by
\[
\rho(T)
\left(
\sum_i\widetilde\pi(a_i)W\eta_i
\right)
\coloneqq
\sum_i\widetilde\pi(a_i)WT\eta_i,
\qquad T\in\mathcal C.
\]
Then
\[
\rho(T)W=WT,
\qquad T\in\mathcal C,
\]
and
\[
\rho(\mathcal C)\subseteq\widetilde\pi(\cL)'.
\]

Since
\[
A',B'\in\pi(\cM)'
\subseteq(\pi\circ\Theta)(\cL)'=\mathcal C,
\]
we may define
\[
\widetilde A'\coloneqq\rho(A'),
\qquad
\widetilde B'\coloneqq\rho(B').
\]
Since $\rho$ is a unital $*$-representation,
we have
\[
\widetilde A',\widetilde B'
\in\widetilde\pi(\cL)'_+,
\qquad
\widetilde A'+\widetilde B'=I_{\widetilde\cK}.
\]
Furthermore,
\[
\widetilde A'W=WA',
\qquad
\widetilde B'W=WB'.
\]

For every $a\in\cL$, we have
\begin{align*}
\widetilde V^*\widetilde A'
\widetilde\pi(a)\widetilde V
&=
V^*W^*\widetilde A'\widetilde\pi(a)WV
=
V^*A'W^*\widetilde\pi(a)WV\\
&=
V^*A'\pi(\Theta(a))V
=
(\Phi\circ\Theta)(a).
\end{align*}
Similarly,
\[
\widetilde V^*\widetilde B'
\widetilde\pi(a)\widetilde V
=
(\Psi\circ\Theta)(a).
\]
Thus $\widetilde A'$ and $\widetilde B'$ are RN derivatives of
$\Phi\circ\Theta$ and $\Psi\circ\Theta$, respectively, in the
Stinespring representation
$(\widetilde\pi,\widetilde V,\widetilde\cK)$.

Moreover, 
$\widetilde A'$ and $\widetilde B'$ commute
with $WW^*$. 
By Lemma~\ref{lem:isometry},
\[
W^*\phi_f(\widetilde A',\widetilde B')W
=
\phi_f(W^*\widetilde A'W,W^*\widetilde B'W)
=
\phi_f(A',B').
\]

Finally, by Proposition~\ref{prop:non-minimal},
we obtain
\begin{align*}
\widehat S_f(\Phi\circ\Theta\|\Psi\circ\Theta)
&\leq
\widetilde V^*
\phi_f(\widetilde A',\widetilde B')
\widetilde V
=
V^*W^*
\phi_f(\widetilde A',\widetilde B')
WV\\
&=
V^*\phi_f(A',B')V
=
\widehat S_f(\Phi\|\Psi).
\end{align*}
\epf

%%%%%%%%%%%%%%%%%%%%%%%%%%%%%%%%%%%%%%%%%%%%%%%%%%%%%%%%%%%%%%%%%%%%%%%%%%%%%%%%%%%%%%%%%%%%%%%%%%%%%%%%%%%%%%%%%%%%%%%%%%%%%%%%%%%%%%%%%%%%%%%%%%%%%%%%%%%%%%%%%%%%%%%%%%%%%%%%%%%%%%%%%%%%%%%%%%%%%%%%%%%%%%%%%%%%%%%%%%%%%%%%%%%%%%%%%%%%%%%%%%%%%%%%%%%%%%%%%%%%%
\subsection{Postcomposition monotonicity}
%%%%%%%%%%%%%%%%%%%%%%%%%%%%%%%%%%%%%%%%%%%%%%%%%%%%%%%%%%%%%%%%%%%%%%%%%%%%%%%%%%%%%%%%%%%%%%%%%%%%%%%%%%%%%%%%%%%%%%%%%%%%%%%%%%%%%%%%%%%%%%%%%%%%%%%%%%%%%%%%%%%%%%%%%%%%%%%%%%%%%%%%%%%%%%%%%%%%%%%%%%%%%%%%%%%%%%%%%%%%%%%%%%%%%%%%%%%%%%%%%%%%%%%%%%%%%%%%%%%%%

Let $\cR$ be a von Neumann algebra acting on a Hilbert space
$\cH_{\cR}$, and let
$\Theta\in\mathrm{CP}(\cN,\cR)$ be normal. Its preadjoint is the
positive linear map
\[
\Theta_*\colon\cR_*\to\cN_*,
\qquad
\Theta_*(\varphi)=\varphi\circ\Theta.
\]
Thus,
\[
(\Theta_*(\varphi))(b)=\varphi(\Theta(b)),
\qquad
b\in\cN,\quad \varphi\in\cR_*,
\]
and
\[
\Theta_*(\cR_*^+)\subseteq\cN_*^+.
\]

For $m\in\widehat\cN_{\mathrm{lb}}$, define
\[
\widehat\Theta(m)(\varphi)
\coloneqq
m(\Theta_*(\varphi))
=
m(\varphi\circ\Theta),
\qquad
\varphi\in\cR_*^+.
\]
Then
\[
\widehat\Theta(m)\in\widehat\cR_{\mathrm{lb}}.
\]
The map
\[
\widehat\Theta\colon
\widehat\cN_{\mathrm{lb}}
\to
\widehat\cR_{\mathrm{lb}}
\]
is positively homogeneous, additive, and order-preserving. Under
the natural embeddings
\[
\cN_{\mathrm{sa}}\subseteq\widehat\cN_{\mathrm{lb}},
\qquad
\cR_{\mathrm{sa}}\subseteq\widehat\cR_{\mathrm{lb}},
\]
it extends $\Theta$:
\[
\widehat\Theta(b)=\Theta(b),
\qquad b\in\cN_{\mathrm{sa}}.
\]
This is the extension used in and immediately before
\cite[Proposition~7.10]{huw_2022}.

\begin{lemma}
\label{lem:extension-compatibility}
Let $\mathcal Q\subseteq\bB(\cK)$ be a von Neumann algebra.
By Proposition~\ref{prop:affiliated}, we identify $\widehat{\mathcal Q}_{\mathrm{lb}}$ with
the set of elements of $\widehat{\bB(\cK)}_{\mathrm{lb}}$ affiliated
with $\mathcal Q$, via $m\mapsto\widetilde m$.
\begin{enumerate}
\item
Let $V\in\bB(\cH,\cK)$ satisfy $V^*\mathcal QV\subseteq\cN$, and put
$\rho(C)\coloneqq V^*CV$ for $C\in\mathcal Q$. Then
$\widehat\rho(m)=V^*mV$ for $m\in\widehat{\mathcal Q}_{\mathrm{lb}}$.
\item
Let $\Xi\colon\mathcal Q\to\cR$ be a normal CP map with Stinespring
representation $(\sigma,\widetilde V,\widetilde\cK)$, where
$\sigma\colon\mathcal Q\to\bB(\widetilde\cK)$ is a normal
$*$-representation. Then
$\widetilde V^*\widehat\sigma(m)\widetilde V=\widehat\Xi(m)$ for
$m\in\widehat{\mathcal Q}_{\mathrm{lb}}$.
\item
If $\rho\colon\mathcal Q\to\cN$ and $\Theta\colon\cN\to\cR$ are normal
CP maps, then
$\widehat{\Theta\circ\rho}=\widehat\Theta\circ\widehat\rho$.
\item
Let $A'\in\mathcal Q$ with $0\leq A'\leq I_{\cK}$ and spectral
measure $E'$. Let $p\colon[0,1]\to(-\infty,+\infty]$ be a Borel
function bounded from below. Then
$p(A')\in\widehat{\bB(\cK)}_{\mathrm{lb}}$, defined by
\[
p(A')(\omega)\coloneqq\int_{[0,1]}p(t)\,d\omega(E'(t)),
\qquad\omega\in\bB(\cK)_*^+,
\]
is affiliated with $\mathcal Q$. Moreover, for every normal
$*$-representation $\sigma$ of $\mathcal Q$,
\[
\widehat\sigma\bigl(p(A')\bigr)=p(\sigma(A')).
\]
\end{enumerate}
\end{lemma}
 
\begin{proof}
(1) Let $\omega\in\bB(\cH)_*^+$. For $C\in\mathcal Q$,
\[
(V\omega V^*)(C)
=
\omega(V^*CV)
=
\omega|_{\cN}(\rho(C))
=
\rho_*(\omega|_{\cN})(C),
\]
so $(V\omega V^*)|_{\mathcal Q}=\rho_*(\omega|_{\cN})$. Hence, under the
identification above,
\[
(V^*mV)(\omega)
=
m\bigl((V\omega V^*)|_{\mathcal Q}\bigr)
=
m\bigl(\rho_*(\omega|_{\cN})\bigr)
=
\widehat\rho(m)(\omega|_{\cN}).
\]
 
(2) Let $\omega\in\bB(\cH_{\cR})_*^+$. For $C\in\mathcal Q$,
\[
(\widetilde V\omega\widetilde V^*)(\sigma(C))
=
\omega(\widetilde V^*\sigma(C)\widetilde V)
=
\omega|_{\cR}(\Xi(C)),
\]
that is,
$\sigma_*(\widetilde V\omega\widetilde V^*)=\Xi_*(\omega|_{\cR})$, where
$\sigma$ is regarded as a map into $\bB(\widetilde\cK)$.
Therefore,
\[
\bigl(\widetilde V^*\widehat\sigma(m)\widetilde V\bigr)(\omega)
=
\widehat\sigma(m)(\widetilde V\omega\widetilde V^*)
=
m\bigl(\Xi_*(\omega|_{\cR})\bigr)
=
\widehat\Xi(m)(\omega|_{\cR}).
\]
 
(3) This follows from
$(\Theta\circ\rho)_*=\rho_*\circ\Theta_*$.
 
(4) Since $E'(\Delta)\in\mathcal Q$ for every Borel set
$\Delta\subseteq[0,1]$, we have
$u'E'(\Delta)u'^*=E'(\Delta)$ for $u'\in\cU(\mathcal Q')$. Hence
$p(A')(u'\omega u'^*)=p(A')(\omega)$, so $p(A')$ is affiliated with
$\mathcal Q$.
 
Since $\sigma$ is normal, $\sigma\circ E'$ is the
spectral measure of $\sigma(A')$. Consequently, for
$\omega\in\bB(\widetilde\cK)_*^+$,
\[
\widehat\sigma\bigl(p(A')\bigr)(\omega)
=
p(A')(\omega\circ\sigma)
=
\int p(t)\,d\omega\bigl(\sigma(E'(t))\bigr)
=
p(\sigma(A'))(\omega).
\]
\end{proof}

\bthm\label{thm:postcomposition}
Let $\Phi,\Psi\in\mathrm{CP}(\cM,\cN)$, and let
$\Theta\in\mathrm{CP}(\cN,\cR)$ be normal. Then
\[
\widehat S_f(\Theta\circ\Phi\|\Theta\circ\Psi)
\leq
\widehat\Theta
\bigl(\widehat S_f(\Phi\|\Psi)\bigr).
\]
\ethm

\bpf
Let $\Gamma\coloneqq\Phi+\Psi$,
and
$
(\pi,V,\cK;A',B')
$
be a minimal SRN realization of $(\Phi,\Psi)$.
Thus,
\[
A'+B'=I_{\cK}.
\]
Let
$
\pi'\colon\cN'\to\bB(\cK)
$
be the unital $*$-representation in
Lemma~\ref{lem:pi-commutant}, and set
$
\mathcal Q\coloneqq\pi'(\cN')'$.
By Lemma~\ref{lem:pi-commutant},
$
A',B'\in\mathcal Q$.

For $C\in\mathcal Q$ and $b'\in\cN'$, we have
\[
(V^*CV)b'
=
V^*C\pi'(b')V
=
V^*\pi'(b')CV
=
b'V^*CV.
\]
Hence
\[
V^*CV\in(\cN')'=\cN.
\]
It follows that
\[
\rho\colon\mathcal Q\to\cN,
\qquad
\rho(C)\coloneqq V^*CV,
\]
is a normal CP map.
Define
\[
\Xi\coloneqq\Theta\circ\rho
\colon\mathcal Q\to\cR.
\]
Since both $\Theta$ and $\rho$ are normal, $\Xi$ is a normal CP
map. Let
\[
(\sigma,\widetilde V,\widetilde\cK)
\]
be a minimal Stinespring representation of $\Xi$. The normality of
$\Xi$ implies that $\sigma$ is a normal $*$-representation of
$\mathcal Q$.

Since
$
\pi'(\cN')\subseteq\pi(\cM)'$,
we have
\[
\pi(\cM)\subseteq\pi'(\cN')'=\mathcal Q.
\]
We may therefore define a unital $*$-representation
\[
\widetilde\pi\colon\cM\to\bB(\widetilde\cK)
\]
by
\[
\widetilde\pi(a)\coloneqq\sigma(\pi(a)),
\qquad a\in\cM.
\]
For every $a\in\cM$,
\begin{align*}
(\Theta\circ\Gamma)(a)
&=
\Theta(V^*\pi(a)V)
=
\Xi(\pi(a))\\
&=
\widetilde V^*\sigma(\pi(a))\widetilde V
=
\widetilde V^*\widetilde\pi(a)\widetilde V.
\end{align*}
Thus
$
(\widetilde\pi,\widetilde V,\widetilde\cK)
$
is a Stinespring representation of $\Theta\circ\Gamma$ that need not be minimal.

Set
\[
\widetilde A'\coloneqq\sigma(A'),
\qquad
\widetilde B'\coloneqq\sigma(B').
\]
Since $\sigma$ is a unital $*$-representation,
we have
\[
0\leq\widetilde A',\widetilde B'
\leq I_{\widetilde\cK},
\qquad
\widetilde A'+\widetilde B'
=
\sigma(A'+B')
=
I_{\widetilde\cK}.
\]
Moreover, since $A'$ and $B'$ commute with $\pi(\cM)$,
\[
\widetilde A',\widetilde B'
\in\widetilde\pi(\cM)'.
\]

For every $a\in\cM$, we have
\begin{align*}
\widetilde V^*\widetilde A'
\widetilde\pi(a)\widetilde V
&=
\widetilde V^*\sigma(A')\sigma(\pi(a))\widetilde V
=
\Xi(A'\pi(a))\\
&=
\Theta(V^*A'\pi(a)V)
=
(\Theta\circ\Phi)(a).
\end{align*}
Similarly,
\[
\widetilde V^*\widetilde B'
\widetilde\pi(a)\widetilde V
=
(\Theta\circ\Psi)(a).
\]
Thus $\widetilde A'$ and $\widetilde B'$ are RN derivatives of
$\Theta\circ\Phi$ and $\Theta\circ\Psi$, respectively, in the
Stinespring representation
$(\widetilde\pi,\widetilde V,\widetilde\cK)$.

Let
\[
p_f(t)\coloneqq\phi_f(t,1-t),
\qquad 0\leq t\leq1.
\]
The function $p_f$ is a Borel function on $[0,1]$ with values in
$(-\infty,+\infty]$, and it is bounded from below. Since
$B'=I_{\cK}-A'$ and $\sigma(B')=I_{\widetilde\cK}-\sigma(A')$, we have
\[
\phi_f(A',B')=p_f(A'),
\qquad
\phi_f(\sigma(A'),\sigma(B'))=p_f(\sigma(A')).
\]
Since $A'\in\mathcal Q$, Lemma~\ref{lem:extension-compatibility}(4)
shows that $m\coloneqq\phi_f(A',B')$ is affiliated with $\mathcal Q$.
We therefore regard $m$ as an element of
$\widehat{\mathcal Q}_{\mathrm{lb}}$ by
Proposition~\ref{prop:affiliated}. The same lemma also gives
\[
\phi_f(\sigma(A'),\sigma(B'))
=
p_f(\sigma(A'))
=
\widehat\sigma\bigl(p_f(A')\bigr)
=
\widehat\sigma(m).
\]
By Lemma~\ref{lem:extension-compatibility}(1)--(3), applied to
$\rho$, $\Xi=\Theta\circ\rho$ and its Stinespring representation
$(\sigma,\widetilde V,\widetilde\cK)$, we have
\[
\widetilde V^*\widehat\sigma(m)\widetilde V
=
\widehat\Xi(m)
=
\widehat\Theta\bigl(\widehat\rho(m)\bigr)
=
\widehat\Theta(V^*mV).
\]

Therefore, by Proposition~\ref{prop:non-minimal}, we obtain
\begin{align*}
\widehat S_f(\Theta\circ\Phi\|\Theta\circ\Psi)
&\leq
\widetilde V^*
\phi_f(\widetilde A',\widetilde B')
\widetilde V
=
\widetilde V^*
\widehat\sigma(m)
\widetilde V\\
&=
\widehat\Theta(V^*mV)
=
\widehat\Theta
\bigl(V^*\phi_f(A',B')V\bigr)
=
\widehat\Theta
\bigl(\widehat S_f(\Phi\|\Psi)\bigr).
\end{align*}
\epf

\bcor
\label{cor:isomorphism-covariance}
Let
$
\alpha\colon\cN\to\cR
$
be a normal $*$-isomorphism. Then, for every
$\Phi,\Psi\in\mathrm{CP}(\cM,\cN)$,
\[
\widehat S_f(\alpha\circ\Phi\|\alpha\circ\Psi)
=
\widehat\alpha
\bigl(\widehat S_f(\Phi\|\Psi)\bigr).
\]
\ecor

\bpf
Theorem~\ref{thm:postcomposition}, applied to $\alpha$, gives
\[
\widehat S_f(\alpha\circ\Phi\|\alpha\circ\Psi)
\leq
\widehat\alpha
\bigl(\widehat S_f(\Phi\|\Psi)\bigr).
\]
Applying the same theorem to $\alpha^{-1}$ and the pair
$(\alpha\circ\Phi,\alpha\circ\Psi)$, and then applying the
order-preserving map $\widehat\alpha$, which satisfies
$\widehat\alpha\circ\widehat{\alpha^{-1}}=\mathrm{id}$ by
Lemma~\ref{lem:extension-compatibility}(3), yields the reverse
inequality.
\epf

\begin{remark}
By Corollary~\ref{cor:isomorphism-covariance}, $\widehat S_f(\Phi\|\Psi)\in\widehat\cN_{\mathrm{lb}}$
does not depend on the choice of a faithful normal representation of $\cN$.
In what follows, we freely choose a convenient one (e.g.\ a standard form or a GNS representation).
\end{remark}

%%%%%%%%%%%%%%%%%%%%%%%%%%%%%%%%%%%%%%%%%%%%%%%%%%%%%%%%%%%%%%%%%%%%%%%%%%%%%%%%%%%%%%%%%%%%%%%%%%%%%%%%%%%%%%%%%%%%%%%%%%%%%%%%%%%%%%%%%%%%%%%%%%%%%%%%%%%%%%%%%%%%%%%%%%%%%%%%%%%%%%%%%%%%%%%%%%%%%%%%%%%%%%%%%%%%%%%%%%%%%%%%%%%%%%%%%%%%%%%%%%%%%%%%%%%%%%%%%%%%%
%%%%%%%%%%%%%%%%%%%%%%%%%%%%%%%%%%%%%%%%%%%%%%%%%%%%%%%%%%%%%%%%%%%%%%%%%%%%%%%%
\subsection{Joint lower semicontinuity}
%%%%%%%%%%%%%%%%%%%%%%%%%%%%%%%%%%%%%%%%%%%%%%%%%%%%%%%%%%%%%%%%%%%%%%%%%%%%%%%%

We first record a consequence of the variational expression for the
PW functional calculus.

\blem
\label{lem:form-lsc}
Let $\cV$ be a complex vector space, and let
$(\rho_i)_{i\in I}$ and $(\sigma_i)_{i\in I}$ be nets of positive
sesquilinear forms on $\cV$. Suppose that
\[
\rho_i(X,Y)\to \rho(X,Y),
\qquad
\sigma_i(X,Y)\to \sigma(X,Y)
\]
for all $X,Y\in\cV$, where $\rho$ and $\sigma$ are positive sesquilinear
forms on $\cV$. Then
\[
q_f[\rho,\sigma](X)
\leq
\liminf_i q_f[\rho_i,\sigma_i](X),
\qquad X\in\cV.
\]
\elem

\bpf
Fix $X\in\cV$. Let
\[
(h_\gamma,A,B,\cH_\gamma)
\]
be the canonical compatible representation of $(\rho,\sigma)$,
where $\gamma=\rho+\sigma$. Thus,
\[
\rho(Y,Z)
=
\langle Ah_\gamma(Y),h_\gamma(Z)\rangle,
\qquad
\sigma(Y,Z)
=
\langle Bh_\gamma(Y),h_\gamma(Z)\rangle,
\]
and
\[
A+B=I_{\cH_\gamma}.
\]

For each $n\in\mathbb N$, let $\alpha_n,\beta_n\in\mathbb R$ and let
$\nu_n$ be the finite positive measure supported on $[1/n,n]$ defined in
\cite[(7.11)--(7.13)]{huw_2022}. These data depend only on $f$ and $n$
(see the paragraph preceding \cite[Lemma~7.9]{huw_2022}). 
By the variational formula
\cite[Theorem~9.4]{huw_2022}, we have
\begin{align*}
q_f[\rho,\sigma](X)
&=
\sup_n\sup_{(\eta, \zeta)}\Bigl[
\alpha_n
\langle Ah_\gamma(X),h_\gamma(X)\rangle
+
\beta_n
\langle Bh_\gamma(X),h_\gamma(X)\rangle\\
&\quad
-
\int_{[1/n,n]}
\frac{1+t}{t}
\bigl(
\langle A\eta(t),\eta(t)\rangle
+
t\langle B\zeta(t),\zeta(t)\rangle
\bigr)
\,d\nu_n(t)\Bigr],
\end{align*}
where the second supremum is taken over pairs of finite-range
$\cH_\gamma$-valued step functions on $[1/n,n]$
satisfying
\[
\eta(t)+\zeta(t)=h_\gamma(X).
\]

We claim that the second supremum is unchanged if $\eta$ and
$\zeta$ are required to take their values in $h_\gamma(\cV)$.

Indeed, let $\eta$ and $\zeta$ be an admissible pair. Passing to a
common finite partition
\[
[1/n,n]=E_1\cup\cdots\cup E_m,
\]
we may assume that
\[
\eta(t)=\eta_j,
\qquad
\zeta(t)=\zeta_j,
\qquad t\in E_j,
\]
where
\[
\eta_j+\zeta_j=h_\gamma(X).
\]
Since $h_\gamma(\cV)$ is dense in $\cH_\gamma$, for every $j$ we
may choose $Y_j^{(k)}\in\cV$ such that
\[
h_\gamma(Y_j^{(k)})\to\eta_j.
\]
Put
\[
Z_j^{(k)}
\coloneqq
X-Y_j^{(k)}.
\]
Then
\[
h_\gamma(Z_j^{(k)})
=
h_\gamma(X)-h_\gamma(Y_j^{(k)})
\to
h_\gamma(X)-\eta_j
=
\zeta_j.
\]
Define step functions $Y^{(k)}$ and $Z^{(k)}$ by
\[
Y^{(k)}(t)=Y_j^{(k)},
\qquad
Z^{(k)}(t)=Z_j^{(k)},
\qquad t\in E_j.
\]
They satisfy
\[
Y^{(k)}(t)+Z^{(k)}(t)=X.
\]
Thus
\[
\bigl(
h_\gamma\circ Y^{(k)},
h_\gamma\circ Z^{(k)}
\bigr)
\]
is an admissible pair in the above variational expression.
Then the quadratic terms corresponding to $Y_j^{(k)}$ and
$Z_j^{(k)}$ converge to those corresponding to $\eta_j$ and
$\zeta_j$. 
Since $\nu_n$ is finite and both
\[
t\mapsto\frac{1+t}{t}
\quad\text{and}\quad
t\mapsto1+t
\]
are bounded on $[1/n,n]$, the corresponding integrals converge.
This proves the claim.

For $X\in\cV$, let $\mathcal P_n(X)$ be the set of all pairs
$(Y,Z)$ of finite-range $\cV$-valued step functions
on $[1/n,n]$ such that
\[
Y(t)+Z(t)=X,
\qquad t\in[1/n,n].
\]
By the claim, we have
\begin{align}
q_f[\rho,\sigma](X)
&=
\sup_{n\in\mathbb N}
\sup_{(Y,Z)\in\mathcal P_n(X)}
F_{n,Y,Z}(\rho,\sigma),
\label{eq:form-variational}
\end{align}
where
\begin{align*}
F_{n,Y,Z}(\rho,\sigma)
&\coloneqq
\alpha_n \rho(X,X)
+
\beta_n \sigma(X,X)\\
&\quad
-
\int_{[1/n,n]}
\frac{1+t}{t}
\Bigl(
\rho(Y(t),Y(t))
+
t\,\sigma(Z(t),Z(t))
\Bigr)
\,d\nu_n(t).
\end{align*}

For fixed $n$ and $(Y,Z)\in\mathcal P_n(X)$, the functions $Y$
and $Z$ have finite range. Thus,
\[
F_{n,Y,Z}(\rho_i,\sigma_i)
\to
F_{n,Y,Z}(\rho,\sigma).
\]
Applying \eqref{eq:form-variational} also to $(\rho_i,\sigma_i)$, we have
\[
F_{n,Y,Z}(\rho_i,\sigma_i)
\leq
q_f[\rho_i,\sigma_i](X).
\]
It follows that
\[
F_{n,Y,Z}(\rho,\sigma)
\leq
\liminf_i q_f[\rho_i,\sigma_i](X).
\]
Taking the supremum first over
$(Y,Z)\in\mathcal P_n(X)$ and then over $n\in\mathbb N$, we obtain
\[
q_f[\rho,\sigma](X)
\leq
\liminf_i q_f[\rho_i,\sigma_i](X).
\]
\epf

\bthm\label{thm:jointconti}
Let $\Phi,\Psi\in\mathrm{CP}(\cM,\cN)$, and let
$(\Phi_i)_{i\in I}$ and $(\Psi_i)_{i\in I}$ be nets in
$\mathrm{CP}(\cM,\cN)$ such that
\[
\Phi_i(a)\to\Phi(a),
\qquad
\Psi_i(a)\to\Psi(a)
\]
pointwise in the $\sigma$-weak topology on $\cN$. 
Then
\[
\widehat S_f(\Phi\|\Psi)(\psi)
\leq
\liminf_i
\widehat S_f(\Phi_i\|\Psi_i)(\psi)
\]
for every $\psi\in\cN_*^+$.
\ethm

\bpf
By Corollary~\ref{cor:isomorphism-covariance}, we may represent
$\cN$ in standard form on $\cH$. Let $\mathcal P\subseteq\cH$ be
the corresponding natural positive cone.

Set
$
\cV\coloneqq\cM\odot\cH$.
For $\Theta\in\mathrm{CP}(\cM,\cN)$, define a positive
sesquilinear form $s_\Theta$ on $\cV$ by
\[
s_\Theta
\left(
\sum_{k=1}^m x_k\otimes\xi_k,
\sum_{l=1}^n y_l\otimes\eta_l
\right)
\coloneqq
\sum_{k=1}^m\sum_{l=1}^n
\left\langle
\Theta(y_l^*x_k)\xi_k,\eta_l
\right\rangle.
\]
The point-$\sigma$-weak convergence of the maps implies that
\[
s_{\Phi_i}(X,Y)\to s_\Phi(X,Y),
\qquad
s_{\Psi_i}(X,Y)\to s_\Psi(X,Y)
\]
for all $X,Y\in\cV$.

Let $\psi\in\cN_*^+$, and let $\xi_\psi\in\mathcal P$ be the vector implementing $\psi$, that is, $\psi=\omega_{\xi_\psi}$.
Set
$
X_\psi\coloneqq1_{\cM}\otimes\xi_\psi$.
Then, by \eqref{eq:form-realization} with $\Gamma=\Phi+\Psi$,
\[
\widehat S_f(\Phi\|\Psi)(\psi)
=
q_f[s_\Phi,s_\Psi](X_\psi),
\]
and similarly,
\[
\widehat S_f(\Phi_i\|\Psi_i)(\psi)
=
q_f[s_{\Phi_i},s_{\Psi_i}](X_\psi).
\]
The assertion follows directly from
Lemma~\ref{lem:form-lsc}.
\epf

\bcor
\label{cor:decreasing-perturbations}
Let $\Phi,\Psi\in\mathrm{CP}(\cM,\cN)$, and let
$(\Theta_n)_{n\geq1}$ be a sequence in
$\mathrm{CP}(\cM,\cN)$ such that
\[
\Theta_{n+1}\leq_{\mathrm{cp}}\Theta_n,
\qquad
\Theta_n(a)\to 0
\]
in the $\sigma$-weak topology for every $a\in\cM$. Then
\[
\widehat S_f(\Phi+\Theta_n\|\Psi+\Theta_n)(\psi)
\to
\widehat S_f(\Phi\|\Psi)(\psi)
\]
for every $\psi\in\cN_*^+$. If $f(1)=0$, the convergence is
increasing.
\ecor

\bpf
First suppose that $f(1)=0$.
Let $\Theta\in\mathrm{CP}(\cM,\cN)$, and let $(\pi,V,\cK)$ be a minimal
Stinespring representation of $\Theta$. Then $\Theta\leq_{\mathrm{cp}}\Theta$,
and the RN derivatives of the pair $(\Theta,\Theta)$ with respect to
$\Theta$ are $A'=B'=I_{\cK}$. Since $\phi_f(1,1)=f(1)$,
Proposition~\ref{prop:independent} gives
\[
\widehat S_f(\Theta\|\Theta)
=
V^*\phi_f(I_{\cK},I_{\cK})V
=
f(1)V^*V
=
f(1)\Theta(1_{\cM})
=
0.
\]
Since $\Theta_{n+1}\leq_{\mathrm{cp}}\Theta_n$, the map
$\Theta_n-\Theta_{n+1}$ belongs to $\mathrm{CP}(\cM,\cN)$, and
\[
\Phi+\Theta_n=(\Phi+\Theta_{n+1})+(\Theta_n-\Theta_{n+1}),
\qquad
\Psi+\Theta_n=(\Psi+\Theta_{n+1})+(\Theta_n-\Theta_{n+1}).
\]
Hence joint subadditivity (Theorem~\ref{thm:subadditivity}) gives
\begin{align*}
\widehat S_f(\Phi+\Theta_n\|\Psi+\Theta_n)
&\leq
\widehat S_f(\Phi+\Theta_{n+1}\|\Psi+\Theta_{n+1})
+
\widehat S_f(\Theta_n-\Theta_{n+1}\|\Theta_n-\Theta_{n+1})\\
&=
\widehat S_f(\Phi+\Theta_{n+1}\|\Psi+\Theta_{n+1}).
\end{align*}
Similarly, 
\[
\widehat S_f(\Phi+\Theta_n\|\Psi+\Theta_n)
\leq
\widehat S_f(\Phi\|\Psi)+\widehat S_f(\Theta_n\|\Theta_n)
=
\widehat S_f(\Phi\|\Psi).
\]
Moreover, $(\Phi+\Theta_n)(a)\to\Phi(a)$ and
$(\Psi+\Theta_n)(a)\to\Psi(a)$ $\sigma$-weakly for every $a\in\cM$.
Combining the above inequalities with Theorem~\ref{thm:jointconti},
we obtain
\[
\widehat S_f(\Phi+\Theta_n\|\Psi+\Theta_n)(\psi)
\uparrow
\widehat S_f(\Phi\|\Psi)(\psi).
\]
 
For a general $f$, set $g\coloneqq f-f(1)$. Then
$g\in\mathrm{OC}(0,+\infty)$ and $g(1)=0$.
Moreover, $g(+0)=f(+0)-f(1)$ and $g'(+\infty)=f'(+\infty)$.
Hence, 
\[
\phi_f(s,t)=\phi_g(s,t)+f(1)\,t,
\qquad (s,t)\in\mathbb R_+^2.
\]
Let $(\pi,V,\cK;A',B')$ be a minimal SRN realization of
$(\Phi,\Psi)$. Since $A'+B'=I_{\cK}$, the spectral calculus of $A'$
gives
\[
\phi_f(A',B')=\phi_g(A',B')+f(1)B'.
\]
Therefore,
\[
\widehat S_f(\Phi\|\Psi)
=
\widehat S_g(\Phi\|\Psi)+f(1)\Psi(1_{\cM}),
\]
and the same identity holds for every pair of CP maps, in particular
for $(\Phi+\Theta_n,\Psi+\Theta_n)$.
Applying the preceding case to $g$ and using
\[
\psi\bigl((\Psi+\Theta_n)(1_{\cM})\bigr)
=
\psi(\Psi(1_{\cM}))+\psi(\Theta_n(1_{\cM}))
\to
\psi(\Psi(1_{\cM})),
\]
we obtain the assertion.
\epf

%%%%%%%%%%%%%%%%%%%%%%%%%%%%%%%%%%%%%%%%%%%%%%%%%%%%%%%%%%%%%%%%%%%%%%%%%%%%%%%%%%%%%%%%%%%%%%%%%%%%%%%%%%%%%%%%%%%%%%%%%%%%%%%%%%%%%%%%%%%%%%%%%%%%%%%%%%%%%%%%%%%%%%%%%%%%%%%%%%%%%%%%%%%%%%%%%%%%%%%%%%%%%%%%%%%%%%%%%%%%%%%%%%%%%%%%%%%%%%%%%%%%%%%%%%%%%%%%%%%%%
\subsection{Martingale convergence}
%%%%%%%%%%%%%%%%%%%%%%%%%%%%%%%%%%%%%%%%%%%%%%%%%%%%%%%%%%%%%%%%%%%%%%%%%%%%%%%%%%%%%%%%%%%%%%%%%%%%%%%%%%%%%%%%%%%%%%%%%%%%%%%%%%%%%%%%%%%%%%%%%%%%%%%%%%%%%%%%%%%%%%%%%%%%%%%%%%%%%%%%%%%%%%%%%%%%%%%%%%%%%%%%%%%%%%%%%%%%%%%%%%%%%%%%%%%%%%%%%%%%%%%%%%%%%%%%%%%%%

The following result extends the martingale convergence theorem for
maximal $f$-divergences of normal positive functionals
\cite[Theorem~5.6]{hiai_2019}.

\bthm\label{thm:martingale}
Let $(\cM_i)_{i\in I}$ be an increasing net of unital
von Neumann subalgebras of $\cM$ such that
\[
\left(\bigcup_{i\in I}\cM_i\right)''=\cM.
\]
Let $\Phi,\Psi\in\mathrm{CP}(\cM,\cN)$ be normal. Then
for every
$\psi\in\cN_*^+$,
\[
\widehat S_f
\bigl(
\Phi|_{\cM_i}\|\Psi|_{\cM_i}
\bigr)(\psi)
\uparrow
\widehat S_f(\Phi\|\Psi)(\psi).
\]
\ethm

\bpf
By Corollary~\ref{cor:isomorphism-covariance}, we may represent
$\cN$ in standard form on $\cH$.

For $i\leq j$, precomposition monotonicity (Theorem~\ref{thm:precomposition}), applied
to the inclusions
\[
\cM_i\hookrightarrow\cM_j\hookrightarrow\cM,
\]
gives
\[
\widehat S_f
\bigl(
\Phi|_{\cM_i}\|\Psi|_{\cM_i}
\bigr)
\leq
\widehat S_f
\bigl(
\Phi|_{\cM_j}\|\Psi|_{\cM_j}
\bigr)
\leq
\widehat S_f(\Phi\|\Psi).
\]
Thus, it remains to prove the reverse inequality in the limit.

Let
$
(\pi,V,\cK;A',B')
$
be a minimal SRN realization of $(\Phi,\Psi)$. Since
$\Phi+\Psi$ is normal, $\pi$ is a normal $*$-representation. For each
$i\in I$, set
\[
\cK_i
\coloneqq
\overline{\pi(\cM_i)V\cH}.
\]
Let $P_i$ be the projection onto $\cK_i$.
Since $(\cM_i)_{i\in I}$ is increasing, so is $(P_i)_{i\in I}$.
Moreover,
since $\bigcup_{i\in I}\cM_i$ is $\sigma$-weakly dense in $\cM$ and $\pi$ is normal,
$P_i\uparrow I_{\cK}$ strongly.

The subspace $\cK_i$ reduces $\pi(\cM_i)$, and
$P_i V=V$. Let
\[
\pi_i(x)
\coloneqq
\pi(x)|_{\cK_i},
\qquad x\in\cM_i.
\]
Then $(\pi_i,V,\cK_i)$ is a minimal Stinespring
representation of $(\Phi+\Psi)|_{\cM_i}$. Its RN derivatives
corresponding to $\Phi|_{\cM_i}$ and
$\Psi|_{\cM_i}$ are
\[
A_i'
\coloneqq
(P_i A'P_i)|_{\cK_i},
\qquad
B_i'
\coloneqq
(P_i B'P_i)|_{\cK_i}.
\]
Consequently,
\[
\widehat S_f
\bigl(
\Phi|_{\cM_i}\|\Psi|_{\cM_i}
\bigr)
=
V^*\phi_f(A_i',B_i')V.
\]

By the direct sum property of the PW functional calculus
\cite[Proposition~4.6]{huw_2022} and $P_i V=V$, we have
\[
\widehat S_f
\bigl(
\Phi|_{\cM_i}\|\Psi|_{\cM_i}
\bigr)(\omega_\xi)
=
\phi_f(P_i A'P_i,
       P_i B'P_i)(\omega_{V\xi})
\]
for every $\xi\in\cH$.

Since $P_i\to I_{\cK}$ strongly, we have
\[
P_i A'P_i\to A',
\qquad
P_i B'P_i\to B'
\]
strongly. Lemma~\ref{lem:form-lsc}, applied to the positive
sesquilinear forms induced by these operators, gives
\begin{align*}
\widehat S_f(\Phi\|\Psi)(\omega_\xi)
&=
\phi_f(A',B')(\omega_{V\xi})\\
&\leq
\liminf_i
\widehat S_f
\bigl(
\Phi|_{\cM_i}\|\Psi|_{\cM_i}
\bigr)(\omega_\xi)
\end{align*}
for every $\xi\in\cH$.

Since $\cN$ is in standard form, every $\psi\in\cN_*^+$ is of the
form $\psi=\omega_\xi|_{\cN}$ for some $\xi\in\cH$. Combining the
last inequality with the monotonicity established at the beginning
of the proof therefore yields the asserted increasing convergence
for every $\psi\in\cN_*^+$.
\epf

%%%%%%%%%%%%%%%%%%%%%%%%%%%%%%%%%%%%%%%%%%%%%%%%%%%%%%%%%%%%%%%%%%%%%%%%%%%%%%%%%%%%%%%
%%%%%%%%%%%%%%%%%%%%%%%%%%%%%%%%%%%%%%%%%%%%%%%%%%%%%%%%%%%%%%%%%%%%%%%%%%%%%%%%%%%%%%%
\subsection{Dense domains and boundedness}
%%%%%%%%%%%%%%%%%%%%%%%%%%%%%%%%%%%%%%%%%%%%%%%%%%%%%%%%%%%%%%%%%%%%%%%%%%%%%%%%%%%%%%%
%%%%%%%%%%%%%%%%%%%%%%%%%%%%%%%%%%%%%%%%%%%%%%%%%%%%%%%%%%%%%%%%%%%%%%%%%%%%%%%%%%%%%%%

We investigate when
$\widehat S_f(\Phi\|\Psi)\in\widehat\cN_{\mathrm{lb}}$
has a dense domain or is represented by a bounded self-adjoint
operator in $\cN$.

We first recall the notions of absolute continuity and singularity
for CP maps from \cite[Definition~7.1]{oka}.

\bdf
Let $\Phi,\Psi\in\mathrm{CP}(\cM,\cN)$.
\begin{enumerate}
\item
We write $\Psi\ll\Phi$ if there exist
$\Psi_n\in\mathrm{CP}(\cM,\cN)$ and $c_n>0$ such that
\[
\Psi_n\uparrow\Psi
\quad\text{in the point-$\sigma$-weak topology},
\qquad
\Psi_n\leq_{\mathrm{cp}}c_n\Phi.
\]

\item
We write $\Phi\perp\Psi$ if the zero map is the only
$\Theta\in\mathrm{CP}(\cM,\cN)$ satisfying
\[
\Theta\leq_{\mathrm{cp}}\Phi,
\qquad
\Theta\leq_{\mathrm{cp}}\Psi.
\]
\end{enumerate}
\edf

Let
$
(\pi,V,\cK;A',B')
$
be a minimal SRN realization of $(\Phi,\Psi)$, and let
$P_\Phi'$ and $P_\Psi'$ be the support projections of $A'$ and
$B'$, respectively. By \cite[Theorem~7.5]{oka},
\[
\Psi\ll\Phi
\iff
P_\Psi'\leq P_\Phi',
\qquad
\Phi\perp\Psi
\iff
P_\Phi'P_\Psi'=0.
\]

Let
\[
A'=\int_{[0,1]}t\,dE'(t)
\]
be the spectral decomposition of $A'$. Since $B'=I_{\cK}-A'$,
\[
P_\Phi'=E'((0,1]),
\qquad
P_\Psi'=E'([0,1)).
\]
Set
\[
e_0'\coloneqq E'(\{0\}),
\qquad
e_1'\coloneqq E'(\{1\}).
\]
It follows that
\[
\Psi\ll\Phi\iff e_0'=0,
\qquad
\Phi\ll\Psi\iff e_1'=0,
\]
and
\[
\Phi\perp\Psi
\iff
e_0'+e_1'=I_{\cK}.
\]

For convenience, set
\[
p_f(t)\coloneqq\phi_f(t,1-t),
\qquad 0\leq t\leq1,
\]
and write $x\vee0\coloneqq\max\{x,0\}$.

The dense-domain assertions below are CP-map analogues of
\cite[Proposition~8.1]{huw_2022}, while the boundedness assertions
are related to \cite[Proposition~8.9]{huw_2022}.

\bprop\label{prop:dense-bounded}
Let $f\in\mathrm{OC}(0,+\infty)$ and
$\Phi,\Psi\in\mathrm{CP}(\cM,\cN)$.

\begin{enumerate}
\item
Suppose that
\[
f'(+\infty)=+\infty>f(+0),
\]
and set
\[
Q_1'
\coloneqq
\int_{[0,1)}
\bigl(p_f(t)\vee0\bigr)\,dE'(t).
\]
Then $\widehat S_f(\Phi\|\Psi)$ has a dense domain if and only if
\[
\Phi\ll\Psi
\quad\text{and}\quad
\operatorname{dom}(Q_1'^{1/2}V)
\ \text{is dense in }\cH.
\]
If $\Phi\leq_{\mathrm{cp}}c\Psi$ for some $c>0$, then
$\widehat S_f(\Phi\|\Psi)$ is bounded.

\item
Suppose that
\[
f'(+\infty)<+\infty=f(+0),
\]
and set
\[
Q_0'
\coloneqq
\int_{(0,1]}
\bigl(p_f(t)\vee0\bigr)\,dE'(t).
\]
Then $\widehat S_f(\Phi\|\Psi)$ has a dense domain if and only if
\[
\Psi\ll\Phi
\quad\text{and}\quad
\operatorname{dom}(Q_0'^{1/2}V)
\ \text{is dense in }\cH.
\]
If $\Psi\leq_{\mathrm{cp}}c\Phi$ for some $c>0$, then
$\widehat S_f(\Phi\|\Psi)$ is bounded.

\item
Suppose that
\[
f'(+\infty)=f(+0)=+\infty,
\]
and set
\[
Q_{0,1}'
\coloneqq
\int_{(0,1)}
\bigl(p_f(t)\vee0\bigr)\,dE'(t).
\]
Then $\widehat S_f(\Phi\|\Psi)$ has a dense domain if and only if
\[
\Phi\ll\Psi,\quad \Psi\ll\Phi
\quad
\text{and}
\quad
\operatorname{dom}(Q_{0,1}'^{1/2}V)
\ \text{is dense in }\cH.
\]
If
$\Phi\leq_{\mathrm{cp}}c\Psi$, 
$\Psi\leq_{\mathrm{cp}}d\Phi$
for some $c,d>0$, then
$\widehat S_f(\Phi\|\Psi)$ is bounded.

\item
If $f'(+\infty)<+\infty$ and $f(+0)<+\infty$,
then $\widehat S_f(\Phi\|\Psi)$ is bounded.
\end{enumerate}
\eprop

\bpf
Since $A'+B'=I_{\cK}$, the spectral functional calculus gives
\[
\phi_f(A',B')
=
p_f(A')
=
\int_{[0,1]}p_f(t)\,dE'(t).
\]

We prove (1). Since $p_f$ is bounded from below on $[0,1]$, for
any $\xi\in\cH$ we have
\[
\widehat S_f(\Phi\|\Psi)(\omega_\xi)<+\infty
\]
if and only if
\[
\xi\in
\ker(e_1'V)
\cap
\operatorname{dom}(Q_1'^{1/2}V).
\]

Consequently, $\widehat S_f(\Phi\|\Psi)$ has a dense domain if
and only if
\[
e_1'V=0
\quad\text{and}\quad
\operatorname{dom}(Q_1'^{1/2}V)
\ \text{is dense in }\cH.
\]
Since $e_1'\in\pi(\cM)'$ and
$\overline{\pi(\cM)V\cH}=\cK$, we have
\[
e_1'V=0\iff e_1'=0.
\]
The support characterization above gives
\[
e_1'=0\iff\Phi\ll\Psi,
\]
which proves the dense domain assertion in (1).

If $\Phi\leq_{\mathrm{cp}}c\Psi$, then the RN order
correspondence gives
\[
A'\leq cB'=c(I_{\cK}-A'),
\]
and hence
\[
A'\leq\frac{c}{1+c}I_{\cK}.
\]
The function $p_f$ is finite and continuous on
\[
\left[0,\frac{c}{1+c}\right],
\]
so $p_f(A')$ is bounded. Therefore,
\[
\widehat S_f(\Phi\|\Psi)=V^*p_f(A')V
\]
is bounded.

(2) follows from (1) by applying it to the transpose
$\widetilde f$ and the ordered pair $(\Psi,\Phi)$, using Proposition~\ref{prop:transpose}.

For (3), the same spectral argument gives
\[
\{\xi\in\cH:
\widehat S_f(\Phi\|\Psi)(\omega_\xi)<+\infty\}
=
\ker(e_0'V)
\cap
\ker(e_1'V)
\cap
\operatorname{dom}(Q_{0,1}'^{1/2}V).
\]
Thus the domain is dense if and only if
$
e_0'=e_1'=0
$
and $\operatorname{dom}(Q_{0,1}'^{1/2}V)$ is dense. These two
endpoint conditions are equivalent to
\[
\Psi\ll\Phi,
\qquad
\Phi\ll\Psi.
\]

If, in addition,
\[
\Phi\leq_{\mathrm{cp}}c\Psi,
\qquad
\Psi\leq_{\mathrm{cp}}d\Phi,
\]
then
\[
\frac{1}{1+d}I_{\cK}
\leq
A'
\leq
\frac{c}{1+c}I_{\cK}.
\]
Since the spectrum of $A'$ is contained in a compact subinterval of
$(0,1)$ on which $p_f$ is bounded, we obtain the boundedness
assertion in (3).

Finally, under the assumptions of (4), $p_f$ is finite and
continuous on $[0,1]$. 
Therefore, $p_f(A')$ is bounded, and so is 
\[
\widehat S_f(\Phi\|\Psi)=V^*p_f(A')V.
\]
\epf

%%%%%%%%%%%%%%%%%%%%%%%%%%%%%%%%%%%%%%%%%%%%%%%%%%%%%%%%%%%%%%%%%%%%%%%%%%%%%%%%%%%%%%%%%%%%%%%%%%%%%%%%%%%%%%%%%%%%%%%%%%%%%%%%%%%%%%%%%%%%%%%%%%%%%%%%%%%%%%%%%%%%%%%%%%%%%%%%%%%%%%%%%%%%%%%%%%%%%%%%%%%%%%%%%%%%%%%%%%%%%%%%%%%%%%%%%%%%%%%%%%%%%%%%%%%%%%%%%%%%%
\section{\texorpdfstring{Agreement with Hiai's maximal $f$-divergence}{Agreement with Hiai's maximal f-divergence}}
%%%%%%%%%%%%%%%%%%%%%%%%%%%%%%%%%%%%%%%%%%%%%%%%%%%%%%%%%%%%%%%%%%%%%%%%%%%%%%%%%%%%%%%%%%%%%%%%%%%%%%%%%%%%%%%%%%%%%%%%%%%%%%%%%%%%%%%%%%%%%%%%%%%%%%%%%%%%%%%%%%%%%%%%%%%%%%%%%%%%%%%%%%%%%%%%%%%%%%%%%%%%%%%%%%%%%%%%%%%%%%%%%%%%%%%%%%%%%%%%%%%%%%%%%%%%%%%%%%%%%

Throughout this section, let
$f\in\mathrm{OC}(0,+\infty)$.

We compare our definition with the maximal $f$-divergence
introduced by Hiai for normal positive functionals in
\cite{hiai_2019}. We denote Hiai's divergence by
\[
\widehat S_f^{\mathrm{Hiai}}(\varphi\|\psi).
\]

Let $\varphi,\psi\in\cM_*^+$ and set
\[
\gamma\coloneqq\varphi+\psi.
\]
Write $D_\omega$ for the Haagerup $L^1$-density of
$\omega\in\cM_*^+$, and let $e'$ be the projection onto
\[
\overline{\cM D_\gamma^{1/2}}\subseteq L^2(\cM).
\]
Since
\[
0\leq\varphi,\psi\leq\gamma,
\]
both $\varphi$ and $\psi$ are strongly absolutely continuous with
respect to $\gamma$. Hence, by \cite[Lemma~3.1]{hiai_2019}, the
operators
\[
A\coloneqq T_{\varphi/\gamma},
\qquad
B\coloneqq T_{\psi/\gamma}
\]
are contractions in $(e'\cM'e')_+$ and satisfy
\[
\varphi(x)
=
\langle AxD_\gamma^{1/2},D_\gamma^{1/2}\rangle,
\qquad
\psi(x)
=
\langle BxD_\gamma^{1/2},D_\gamma^{1/2}\rangle,
\qquad x\in\cM.
\]
By the uniqueness assertion in \cite[Lemma~3.1]{hiai_2019},
\[
A+B=e'.
\]
Thus, regarded as operators on $e'L^2(\cM)$,
\[
A+B=I_{e'L^2(\cM)}.
\]

Since $B=I_{e'L^2(\cM)}-A$, the operators $A$ and $B$ commute.
Let
\[
A=\int_{[0,1]}t\,dE(t)
\]
be the spectral decomposition of $A$ on $e'L^2(\cM)$. By the
commuting-pair property of the PW functional calculus
\cite[Definition~4.1(1) and Theorem~4.3]{huw_2022},
\begin{align*}
\phi_f(A,B)(\omega_{D_\gamma^{1/2}})
&=
\int_{[0,1]}
\phi_f(t,1-t)\,
d\|E(t)D_\gamma^{1/2}\|^2.
\end{align*}
By \cite[Theorem~4.2]{hiai_2019}, the right-hand side is equal to
Hiai's maximal $f$-divergence. Therefore,
\[
\widehat S_f^{\mathrm{Hiai}}(\varphi\|\psi)
=
\phi_f(A,B)(\omega_{D_\gamma^{1/2}}).
\]

\bthm\label{thm:hiai}
Let $f\in\mathrm{OC}(0,+\infty)$ and
$\varphi,\psi\in\cM_*^+$, and define normal CP maps
$\Phi,\Psi\colon\cM\to\mathbb C$ by
\[
\Phi(x)\coloneqq\varphi(x),
\qquad
\Psi(x)\coloneqq\psi(x).
\]
Under the canonical identification
$
\widehat{\mathbb C}_{\mathrm{lb}}
\simeq(-\infty,+\infty]$,
we have
\[
\widehat S_f(\Phi\|\Psi)
=
\widehat S_f^{\mathrm{Hiai}}(\varphi\|\psi).
\]
\ethm

\bpf
Set $\gamma\coloneqq\varphi+\psi$. If $\gamma=0$, then
$\varphi=\psi=0$, and both sides of the asserted equality are zero.
Hence, assume that $\gamma\neq0$.

Let $(\cH_\gamma,\pi_\gamma,\xi_\gamma)$ be the GNS
representation of $\gamma$, and define
\[
V\colon\mathbb C\to\cH_\gamma,
\qquad
Vz\coloneqq z\xi_\gamma.
\]
Then $(\pi_\gamma,V,\cH_\gamma)$ is a minimal Stinespring
representation of $\Phi+\Psi$. Let
$A',B'\in\pi_\gamma(\cM)'$ be the corresponding RN derivatives.

Define
\[
U_0\colon\pi_\gamma(\cM)\xi_\gamma\to e'L^2(\cM),
\qquad
U_0(\pi_\gamma(x)\xi_\gamma)
\coloneqq
xD_\gamma^{1/2}.
\]
For $x,y\in\cM$,
\[
\langle\pi_\gamma(x)\xi_\gamma,\pi_\gamma(y)\xi_\gamma\rangle
=
\gamma(y^*x)
=
\langle xD_\gamma^{1/2},yD_\gamma^{1/2}\rangle.
\]
Thus $U_0$ is well defined and isometric. Since its range is dense
in $e'L^2(\cM)$, it extends to a unitary
\[
U\colon\cH_\gamma\to e'L^2(\cM)
\]
satisfying
\[
U\xi_\gamma=D_\gamma^{1/2}.
\]

Moreover, for $x,y\in\cM$,
\begin{align*}
\langle UA'U^*xD_\gamma^{1/2},yD_\gamma^{1/2}\rangle
&=
\langle A'\pi_\gamma(x)\xi_\gamma,
        \pi_\gamma(y)\xi_\gamma\rangle\\
&=
\varphi(y^*x)\\
&=
\langle AxD_\gamma^{1/2},yD_\gamma^{1/2}\rangle.
\end{align*}
Hence
$
UA'U^*=A$.
Similarly,
$
UB'U^*=B$.

By unitary covariance of the PW functional calculus,
\[
\widehat S_f(\Phi\|\Psi)
=
\phi_f(A',B')(\omega_{\xi_\gamma})
=
\phi_f(A,B)(\omega_{D_\gamma^{1/2}})
=
\widehat S_f^{\mathrm{Hiai}}(\varphi\|\psi).
\]
\epf

\brem
Hiai proved that
$\widehat S_f^{\mathrm{Hiai}}$ is jointly lower semicontinuous in
the norm topology
\cite[Theorem~5.5]{hiai_2019}. He further observed that joint lower
semicontinuity in the $\sigma(\cM_*,\cM)$-topology follows when
$\cM$ is injective, while the case of a general von Neumann algebra
was left open; see \cite[Remark~5.7]{hiai_2019}.

For the BS case, Hollands and Ranallo later
proved lower semicontinuity
using their variational formula; see
\cite[Proposition~2.17 and Remark~2.18]{h-r_2023}. 

Under the identification of normal positive functionals with normal
CP maps from $\cM$ to $\mathbb C$,
Theorems~\ref{thm:jointconti} and~\ref{thm:hiai} imply that
\[
(\varphi,\psi)
\mapsto
\widehat S_f^{\mathrm{Hiai}}(\varphi\|\psi)
\]
is jointly lower semicontinuous with respect to the product
$\sigma(\cM_*,\cM)$-topology.
\erem

%%%%%%%%%%%%%%%%%%%%%%%%%%%%%%%%%%%%%%%%%%%%%%%%%%%%%%%%%%
%%%%%%%%%%%%%%%%%%%%%%%%%%%%%%%%%%%%%%%%%%%%%%%%%%%%%%%%%%
%%%%%%%%%%%%%%%%%%%%%%%%%%%%%%%%%%%%%%%%%%%%%%%%%%%%%%%%%%
%%%%%%%%%%%%%%%%%%%%%%%%%%%%%%%%%%%%%%%%%%%%%%%%%%%%%%%%%%
\section{Special cases and finite-dimensional examples}
%%%%%%%%%%%%%%%%%%%%%%%%%%%%%%%%%%%%%%%%%%%%%%%%%%%%%%%%%%
%%%%%%%%%%%%%%%%%%%%%%%%%%%%%%%%%%%%%%%%%%%%%%%%%%%%%%%%%%
%%%%%%%%%%%%%%%%%%%%%%%%%%%%%%%%%%%%%%%%%%%%%%%%%%%%%%%%%%
%%%%%%%%%%%%%%%%%%%%%%%%%%%%%%%%%%%%%%%%%%%%%%%%%%%%%%%%%%

%%%%%%%%%%%%%%%%%%%%%%%%%%%%%%%%%%%%%%%%%%%%%%%%%%%%%%%%%%%%%%%%%%%%%%%%%%%%%%%%%%%%%%%%%%%%%%%%%%%%%%%%%%%%%%%%%%%%%%%%%%%%%%%%%%%%%%%%%%%%%%%%%
%%%%%%%%%%%%%%%%%%%%%%%%%%%%%%%%%%%%%%%%%%%%%%%%%%%%%%%%%%
\subsection{Belavkin--Staszewski relative entropy}
%%%%%%%%%%%%%%%%%%%%%%%%%%%%%%%%%%%%%%%%%%%%%%%%%%%%%%%%%%%%%%%%%%%%%%%%%%%%%%%%%%%%%%%%%%%%%%%%%%%%%%%%%%%%%%%%%%%%%%%%%%%%%%%%%%%%%%%%%%%%%%%%%
%%%%%%%%%%%%%%%%%%%%%%%%%%%%%%%%%%%%%%%%%%%%%%%%%%%%%%%%%%

In this subsection,
let $\eta(t)\coloneqq t\log t$. For $A, B\in\bB(\cK)_{++}$,
the corresponding perspective is the negative of the relative
operator entropy:
\[
\phi_\eta(A,B)
=
-S(A\mid B),
\]
where
\[
S(A\mid B)
\coloneqq
A^{1/2}
\log\left(A^{-1/2}BA^{-1/2}\right)
A^{1/2}.
\]
Moreover, we have
\[
-S(A\mid B)
=
\lim_{\alpha\downarrow0}
\frac{A-A\#_\alpha B}{\alpha};
\]
see \cite{fk}.

The PW functional calculus extends this formula to arbitrary
positive operators. More precisely, by
\cite[Example~8.12]{huw_2022},
\[
\phi_\eta(A,B)(\omega)
=
\lim_{\alpha\downarrow0}
\omega\left(
\frac{A-A\#_\alpha B}{\alpha}
\right)
\]
increasingly for every $\omega\in\bB(\cK)_*^+$.

\begin{definition}
The \emph{operator-valued Belavkin--Staszewski relative entropy}
of $\Phi$ with respect to $\Psi$ is defined by
\[
D_{\mathrm{BS}}^{\mathrm{op}}(\Phi\|\Psi)
\coloneqq
\widehat S_\eta(\Phi\|\Psi)
\in\widehat\cN_{\mathrm{lb}}.
\]
\end{definition}

Fix a minimal SRN realization
$(\pi,V,\cK;A',B')$ of $(\Phi,\Psi)$.
For $0<\alpha<1$, let $\Phi\#_\alpha\Psi$ be the weighted
geometric mean of CP maps defined by
\[
(\Phi\#_\alpha\Psi)(x)
=
V^*(A'\#_\alpha B')\pi(x)V,
\qquad x\in\cM;
\]
see \cite[Example~6.8]{oka}.
Then, for every $\psi\in\cN_*^+$,
\[
D_{\mathrm{BS}}^{\mathrm{op}}(\Phi\|\Psi)(\psi)
=
\lim_{\alpha\downarrow0}
\psi\left(
\frac{
\Phi(1_{\cM})-(\Phi\#_\alpha\Psi)(1_{\cM})
}{\alpha}\right)
\]
increasingly. Indeed, we may choose $\omega\in\bB(\cH)_*^+$ with
$\omega|_{\cN}=\psi$. By Proposition~\ref{prop:affiliated} and
\cite[Example~8.12]{huw_2022} applied to
$V\omega V^*\in\bB(\cK)_*^+$, we have
\begin{align*}
D_{\mathrm{BS}}^{\mathrm{op}}(\Phi\|\Psi)(\psi)
&=
\phi_\eta(A',B')(V\omega V^*)
=
\lim_{\alpha\downarrow0}
\omega\left(
V^*\frac{A'-A'\#_\alpha B'}{\alpha}V
\right),
\end{align*}
and
$V^*A'V=\Phi(1_{\cM})$,
$V^*(A'\#_\alpha B')V=(\Phi\#_\alpha\Psi)(1_{\cM})$.

By Theorem~\ref{thm:hiai}, this definition agrees with the usual
BS relative entropy for normal positive
functionals \cite{bs,hiai_2019}. Moreover, if
\[
\Phi\leq_{\mathrm{cp}}c\Psi
\]
for some $c>0$, then
$D_{\mathrm{BS}}^{\mathrm{op}}(\Phi\|\Psi)$ is bounded by
Proposition~\ref{prop:dense-bounded}.

%%%%%%%%%%%%%%%%%%%%%%%%%%%%%%%%%%%%%%%%%%%%%%%%%%%%%%%%%%%%%%%%%%%%%%%%%%%%%%%%%%%%%%%%%%%%%%%%%%%%%%%%%%%%%%%%%%%%%%%%%%%%%%%%%%%%%%%%%%%%%%%%%
%%%%%%%%%%%%%%%%%%%%%%%%%%%%%%%%%%%%%%%%%%%%%%%%%%%%%%%%%%
\subsection{Finite-dimensional examples}
%%%%%%%%%%%%%%%%%%%%%%%%%%%%%%%%%%%%%%%%%%%%%%%%%%%%%%%%%%%%%%%%%%%%%%%%%%%%%%%%%%%%%%%%%%%%%%%%%%%%%%%%%%%%%%%%%%%%%%%%%%%%%%%%%%%%%%%%%%%%%%%%%
%%%%%%%%%%%%%%%%%%%%%%%%%%%%%%%%%%%%%%%%%%%%%%%%%%%%%%%%%%

Throughout this subsection, let
$f\in\mathrm{OC}(0,+\infty)$.

\bex\label{ex:preparation}
Let $m,n\in\mathbb N$, and for $k\in\mathbb N$ let
$\operatorname{tr}_k$ denote the normalized trace on
$M_k(\mathbb C)$. 
Let
\[
\cH=L^2(M_n(\mathbb C),\operatorname{tr}_n),
\]
where $M_n(\mathbb C)$ acts on $\cH$ by left multiplication.

For $\rho,\sigma\in M_n(\mathbb C)_{++}$, define
$
\Phi,\Psi\colon M_m(\mathbb C)\to M_n(\mathbb C)
$
by
\[
\Phi(X)=\operatorname{tr}_m(X)\rho,
\qquad
\Psi(X)=\operatorname{tr}_m(X)\sigma.
\]
Set
$
C\coloneqq\rho+\sigma
$
and
\[
A\coloneqq C^{-1/2}\rho C^{-1/2},
\qquad
B\coloneqq C^{-1/2}\sigma C^{-1/2}.
\]
Then
\[
A, B>0,
\qquad
A+B=I_n.
\]

Let
\[
\cK
\coloneqq
L^2(M_m(\mathbb C),\operatorname{tr}_m)\otimes\cH,
\]
and define
\[
\pi(X)\coloneqq L_X\otimes I_{\cH},
\qquad X\in M_m(\C),
\]
and
\[
V\xi\coloneqq I_m\otimes C^{1/2}\xi,
\qquad \xi\in\cH.
\]
Then
\[
V^*\pi(X)V
=
\operatorname{tr}_m(X)C
=
(\Phi+\Psi)(X).
\]
Since $C$ is invertible, 
the Stinespring representation $(\pi, V, \cK)$ is minimal.

The corresponding RN derivatives are
\[
A'=I\otimes A,
\qquad
B'=I\otimes B.
\]

It follows that
\[
\phi_f(A',B')
=
I\otimes\phi_f(A, B),
\]
and hence
\begin{align*}
\widehat S_f(\Phi\|\Psi)
&=
C^{1/2}
\phi_f(A, B)
C^{1/2}\\
&=
\phi_f(\rho,\sigma),
\end{align*}
where the last equality follows from the operator homogeneity of
the PW functional calculus. Since $\sigma$ is invertible,
\[
\widehat S_f(\Phi\|\Psi)
=
\sigma^{1/2}
f(\sigma^{-1/2}\rho\sigma^{-1/2})
\sigma^{1/2}.
\]
Let $\operatorname{Tr}_n$ be the
non-normalized trace.
Then
\[
\widehat S_f(\Phi\|\Psi)(\operatorname{Tr}_n)
=
\operatorname{Tr}_n
\left(
\sigma^{1/2}
f(\sigma^{-1/2}\rho\sigma^{-1/2})
\sigma^{1/2}
\right),
\]
which is the usual matrix maximal $f$-divergence
\cite{mat,hiai_2019}.
\eex

%%%%%%%%%%%%%%%%%%%%%%%%%%%%%%%%%%%%%%%%%%%%%%%%%%%%%%%%%%%%%%%%%%%%%%%%%%%%%%%%%%%%%%%%%%%%%%%%%%%%%%%%%%%%%%%%%%%%%%%%%%%%%%%%%%%%%%%%%%%%%%%%%%%%%%%%%%%%%%%%%%%%%%%%%%%%%%%%%%%%%%%%%%%%%%%%%%%%%%%%%%%%%%%%%%%%%%%%%%%%%%%%%%%%%%%%

\bex[Schur multipliers]\label{ex:schur}
Let $\cM=\cN=M_n(\C)$, and represent $\cN$ on
\[
\cH=L^2(M_n(\C),\operatorname{Tr}_n)=M_n(\C)
\]
by left multiplication.
Let $\{e_1,\ldots,e_n\}$ be the canonical
orthonormal basis of $\C^n$, and set
\[
E_{ij}=e_ie_j^*,
\qquad 1\leq i,j\leq n.
\]
Then $\{E_{ij}\}_{i,j=1}^n$ is both a system of matrix units and an
orthonormal basis of $\cH$.

Let
\[
A=[a_{ij}],\quad B=[b_{ij}]\in M_n(\C)_{++},
\]
and define the Schur multipliers
\[
\Phi(X)=A\circ X,
\qquad
\Psi(X)=B\circ X,
\qquad X\in M_n(\C),
\]
where $\circ$ denotes the Schur product. 
Recall that $\Phi$ and $\Psi$ are CP maps. Set
$
\Gamma\coloneqq\Phi+\Psi$, and 
$
C\coloneqq A+B.
$
Thus,
\[
\Gamma(X)=C\circ X,
\qquad X\in M_n(\C).
\]

Equip the vector space $\C^n$ with the inner product
\[
\langle\xi,\eta\rangle_{\cH_C}
\coloneqq
\langle C\xi,\eta\rangle_{\C^n},
\]
and denote the resulting Hilbert space by $\cH_C$. Let
\[
T\colon\C^n\to\cH_C,
\qquad
T\xi\coloneqq\xi.
\]
Since $C$ is invertible, the map $T$ is a linear isomorphism, which is not unitary unless $C=I_n$.
Moreover,
\[
\langle T^*Te_j,e_i\rangle_{\C^n}
=
\langle e_j,e_i\rangle_{\cH_C}
=
c_{ij},
\]
and hence
\[
T^*T=C.
\]

Set
\[
\cK=\C^n\otimes\cH_C\otimes\C^n.
\]
Define a unital $*$-representation
$
\pi\colon M_n(\C)\to\bB(\cK)
$
and an operator $V\colon\cH\to\cK$ by
\[
\pi(X)=X\otimes I_{\cH_C}\otimes I_{\C^n},
\qquad
VE_{ij}=e_i\otimes e_i\otimes e_j.
\]

For $X=[x_{ij}]\in M_n(\C)$, we have
\begin{align*}
\left\langle
V^*\pi(X)VE_{kl},E_{ij}
\right\rangle_{\cH}
&=
\left\langle
(Xe_k)\otimes e_k\otimes e_l,
e_i\otimes e_i\otimes e_j
\right\rangle_{\cK}
\\
&=
\langle Xe_k,e_i\rangle_{\C^n}
\langle e_k,e_i\rangle_{\cH_C}
\langle e_l,e_j\rangle_{\C^n}
\\
&=
x_{ik}c_{ik}\delta_{jl}
\\
&=
\left\langle
(C\circ X)E_{kl},E_{ij}
\right\rangle_{\cH}.
\end{align*}
Therefore,
\[
V^*\pi(X)V=\Gamma(X).
\]
Moreover,
\[
\pi(E_{ij})VE_{jk}
=
e_i\otimes e_j\otimes e_k,
\]
and these vectors span $\cK$. Thus $(\pi,V,\cK)$ is a minimal
Stinespring representation of $\Gamma$.

Define $\widehat A,\widehat B\in\bB(\cH_C)$ by
\[
T^*\widehat AT=A,
\qquad
T^*\widehat BT=B.
\]
Equivalently,
\[
\langle\widehat Ae_j,e_i\rangle_{\cH_C}=a_{ij},
\qquad
\langle\widehat Be_j,e_i\rangle_{\cH_C}=b_{ij}.
\]
Since $A,B>0$ and $T$ is invertible, we have
$
\widehat A,\widehat B>0$.
Furthermore,
\[
T^*(\widehat A+\widehat B)T
=
A+B
=
C
=
T^*T,
\]
and hence
\[
\widehat A+\widehat B=I_{\cH_C}.
\]

Set
\[
A'
\coloneqq
I_{\C^n}\otimes\widehat A\otimes I_{\C^n},
\qquad
B'
\coloneqq
I_{\C^n}\otimes\widehat B\otimes I_{\C^n}.
\]
Then
\[
A',B'\in\pi(M_n(\C))',
\qquad
A',B'>0,
\qquad
A'+B'=I_{\cK}.
\]
For $X=[x_{ij}]\in M_n(\C)$, we obtain
\begin{align*}
\left\langle
V^*A'\pi(X)VE_{kl},E_{ij}
\right\rangle_{\cH}
&=
\left\langle
(Xe_k)\otimes(\widehat Ae_k)\otimes e_l,
e_i\otimes e_i\otimes e_j
\right\rangle_{\cK}
\\
&=
x_{ik}a_{ik}\delta_{jl}
\\
&=
\left\langle
(A\circ X)E_{kl},E_{ij}
\right\rangle_{\cH}.
\end{align*}
Therefore,
\[
\Phi(X)=V^*A'\pi(X)V.
\]
Similarly,
\[
\Psi(X)=V^*B'\pi(X)V.
\]
Thus $A'$ and $B'$ are the RN derivatives of $\Phi$ and $\Psi$,
respectively.

It follows that
\[
\phi_f(A',B')
=
I_{\C^n}\otimes
\phi_f(\widehat A,\widehat B)
\otimes I_{\C^n}.
\]
In particular, $\phi_f(A',B')$ is bounded because
$\widehat A$ and $\widehat B$ are positive and invertible.

Let
\[
\xi=\sum_{i,j=1}^n\xi_{ij}E_{ij}\in\cH.
\]
Then
\begin{align*}
\widehat S_f(\Phi\|\Psi)(\omega_\xi)
&=
\left\langle
\phi_f(A',B')V\xi,V\xi
\right\rangle_{\cK}
\\
&=
\sum_{i,j,k,l}
\xi_{ij}\overline{\xi_{kl}}\,
\langle e_i,e_k\rangle_{\C^n}
\left\langle
\phi_f(\widehat A,\widehat B)e_i,e_k
\right\rangle_{\cH_C}
\langle e_j,e_l\rangle_{\C^n}
\\
&=
\sum_{i,j=1}^n
|\xi_{ij}|^2
\left\langle
\phi_f(\widehat A,\widehat B)e_i,e_i
\right\rangle_{\cH_C}.
\end{align*}

Since $T$ is invertible, the transformer equality for the PW
functional calculus gives
\[
T^*\phi_f(\widehat A,\widehat B)T
=
\phi_f(T^*\widehat AT,T^*\widehat BT)
=
\phi_f(A,B).
\]
Since $Te_i=e_i$, it follows
that
\[
\left\langle
\phi_f(\widehat A,\widehat B)e_i,e_i
\right\rangle_{\cH_C}
=
\left\langle
T^*\phi_f(\widehat A,\widehat B)Te_i,e_i
\right\rangle_{\C^n}
=
\left\langle
\phi_f(A,B)e_i,e_i
\right\rangle_{\C^n}.
\]
Consequently,
\[
\widehat S_f(\Phi\|\Psi)(\omega_\xi)
=
\sum_{i,j=1}^n
|\xi_{ij}|^2
\left\langle
\phi_f(A,B)e_i,e_i
\right\rangle_{\C^n}.
\]

Let $\cD_n\subseteq M_n(\C)$ be the diagonal maximal abelian
subalgebra associated with $\{e_1,\ldots,e_n\}$, and let
\[
E_{\cD_n}(X)
\coloneqq
\sum_{i=1}^n
\langle Xe_i,e_i\rangle E_{ii}
\]
be the canonical conditional expectation onto $\cD_n$. Since
$\cN$ acts on $\cH$ by left multiplication, the preceding identity
shows that
\[
\widehat S_f(\Phi\|\Psi)
=
E_{\cD_n}\bigl(\phi_f(A,B)\bigr)
\in\cD_n.
\]

Since $E_{\cD_n}$ preserves the canonical trace, we obtain
\[
\widehat S_f(\Phi\|\Psi)(\operatorname{Tr}_n)
=
\operatorname{Tr}_n
\left(
E_{\cD_n}(\phi_f(A,B))
\right)
=
\operatorname{Tr}_n(\phi_f(A,B)).
\]
Hence the operator-valued maximal $f$-divergence of two Schur multipliers is the
diagonal part of the operator perspective of their symbol matrices,
while its value at the canonical trace recovers their usual matrix
maximal $f$-divergence.
\eex

%%%%%%%%%%%%%%%%%%%%%%%%%%%%%%%%%%%%%%%%%%%%%%%%%%%%%%%%%%%%%%%%%%%%%%%%%%%%%%%%%%%%%%%%%%%%%%%%%%%%%%%%%%%%%%%%%%%%%%%%%%%%%%%%%%%%%%%%%%%%%%%%%%%%%%%%%%%%%%%%%%%%%%%%%%%%%%%%%%%%%%%%%%%%%%%%%%%%%%%%%%%%%%%%%%%%%%%%%%%%%%%%%%%%%%%%

\bex[Congruence maps]\label{ex:congruence}
Let $n\geq 2$ and $\cM=\cN=M_n(\C)$.
Assume that $\cN$ acts on
$\cH=\C^n$ in the standard way. Take
\[
A,B\in M_n(\C)_{++},
\]
and assume that $A$ and $B$ are not proportional, that is,
\[
B\neq\lambda A
\]
for every $\lambda>0$. Define
$\Phi,\Psi\in\mathrm{CP}(M_n(\C),M_n(\C))$ by
\[
\Phi(X)=A^{1/2}XA^{1/2},
\qquad
\Psi(X)=B^{1/2}XB^{1/2},
\qquad X\in M_n(\C),
\]
and set
$
\Gamma\coloneqq\Phi+\Psi$.

Let $\cK=\cH\oplus\cH$, and define
$\pi\colon M_n(\C)\to\bB(\cK)$ and
$V\colon\cH\to\cK$ by
\[
\pi(X)=
\begin{bmatrix}
X&0\\
0&X
\end{bmatrix},
\qquad
V\xi=
\begin{bmatrix}
A^{1/2}\xi\\
B^{1/2}\xi
\end{bmatrix}.
\]
Then
\[
V^*
\begin{bmatrix}
\xi\\
\eta
\end{bmatrix}
=
A^{1/2}\xi+B^{1/2}\eta,
\]
and hence
\[
V^*\pi(X)V
=
A^{1/2}XA^{1/2}
+
B^{1/2}XB^{1/2}
=
\Gamma(X).
\]
Thus $(\pi,V,\cK)$ is a Stinespring representation of $\Gamma$.

We next verify that it is minimal. First note that
$A^{1/2}$ and $B^{1/2}$ are linearly independent. 
Indeed, suppose that $aA^{1/2}+bB^{1/2}=0$ for some $a,b\in\C$.
If $b=0$, then $a=0$, since $A^{1/2}$ is invertible.
If $b\neq0$, then $B^{1/2}=cA^{1/2}$ with $c\coloneqq-a/b\neq0$,
and hence $B=(B^{1/2})^*B^{1/2}=|c|^2A$, a contradiction.

Let $P$ be the projection onto
\[
\overline{\pi(M_n(\C))V\cH}.
\]
This subspace reduces $\pi(M_n(\C))$, and hence
\[
P\in\pi(M_n(\C))'
=
M_2(\C)\otimes I_n.
\]
Suppose that $P\neq I_{\cK}$. Then
\[
Q\coloneqq I_{\cK}-P
=
\begin{bmatrix}
q_{11}I_n&q_{12}I_n\\
q_{21}I_n&q_{22}I_n
\end{bmatrix}
\]
is a nonzero projection satisfying $QV=0$. Consequently,
\[
q_{11}A^{1/2}+q_{12}B^{1/2}=0,
\qquad
q_{21}A^{1/2}+q_{22}B^{1/2}=0.
\]
The linear independence of $A^{1/2}$ and $B^{1/2}$ implies that
$q_{ij}=0$ for all $i,j$, contradicting $Q\neq0$. Thus
\[
P=I_{\cK},
\]
and $(\pi,V,\cK)$ is a minimal Stinespring representation of
$\Gamma$.

Set
\[
A'=
\begin{bmatrix}
I_n&0\\
0&0
\end{bmatrix},
\qquad
B'=
\begin{bmatrix}
0&0\\
0&I_n
\end{bmatrix}.
\]
Then
\[
A',B'\in\pi(M_n(\C))',
\qquad
A'+B'=I_{\cK},
\]
and
\[
V^*A'\pi(X)V
=
A^{1/2}XA^{1/2}
=
\Phi(X),
\qquad
V^*B'\pi(X)V
=
B^{1/2}XB^{1/2}
=
\Psi(X).
\]
Therefore, $A'$ and $B'$ are the RN derivatives of $\Phi$ and
$\Psi$, respectively.

Since $A'$ and $B'$ are mutually orthogonal projections, their
support projections are mutually orthogonal. 
The support characterization of singularity in \cite[Theorem~7.5]{oka} yields
$\Phi\perp\Psi$.
Moreover,
$
A'\wedge B'=0$.
Hence, by 
\cite[Proposition~9.3]{huw_2022},
we have
\[
\phi_f(A',B')
=
f'(+\infty)A'+f(+0)B'
\]
in $\widehat{\bB(\cK)}_{\mathrm{lb}}$. Consequently, for every
$\omega\in M_n(\C)_*^+$,
\begin{align*}
\widehat S_f(\Phi\|\Psi)(\omega)
&=
\bigl(V^*\phi_f(A',B')V\bigr)(\omega)\\
&=
f'(+\infty)\,\omega(V^*A'V)
+
f(+0)\,\omega(V^*B'V)\\
&=
f'(+\infty)\,\omega(A)
+
f(+0)\,\omega(B),
\end{align*}
with the usual convention
$
(+\infty)\cdot0=0
$. 
Equivalently,
\[
\widehat S_f(\Phi\|\Psi)
=
f'(+\infty)A+f(+0)B
\]
in $\widehat{M_n(\C)}_{\mathrm{lb}}$.
In particular,
\[
\widehat S_f(\Phi\|\Psi)(\operatorname{Tr}_n)
=
f'(+\infty)\operatorname{Tr}_n(A)
+
f(+0)\operatorname{Tr}_n(B).
\]
\eex

\brem
Examples~\ref{ex:schur} and~\ref{ex:congruence}
illustrate two different types of overlap between CP maps, as
described by the support projections of their RN derivatives in
a common minimal Stinespring representation.

For the Schur multipliers in Example~\ref{ex:schur}, 
the RN derivatives $A'$ and $B'$ (equivalently, $\widehat A$ and $\widehat B$) are invertible.
Consequently, their support projections are both equal to the
identity, and the two Schur multipliers are mutually absolutely
continuous.
In this case,
\[
\widehat S_f(\Phi\|\Psi)(\operatorname{Tr}_n)
=
\operatorname{Tr}_n\bigl(\phi_f(A,B)\bigr),
\]
so the divergence retains the relative geometry of the two symbol
matrices $A$ and $B$.

By contrast, for the congruence maps in
Example~\ref{ex:congruence}, the RN derivatives are mutually
orthogonal projections whenever $A$ and $B$ are not proportional.
Hence,
\[
\Phi\perp\Psi,
\]
and
\[
\widehat S_f(\Phi\|\Psi)
=
f'(+\infty)A+f(+0)B.
\]
Thus, the divergence is determined entirely by the boundary values
of the perspective.

The additional point in Example~\ref{ex:congruence} is that the
singularity occurs at the level of the CP maps. Both
\[
\Phi(1_n)=A,
\qquad
\Psi(1_n)=B
\]
are invertible, and $A$ and $B$ may be arbitrarily close in norm,
while $\Phi$ and $\Psi$ remain mutually singular whenever $A$ and
$B$ are not proportional.

Therefore, the operator-valued maximal $f$-divergence of CP maps
is governed by the relative position of their RN derivatives in a
common minimal Stinespring representation. In general, this
relative position cannot be determined solely from the values of
the maps at the identity.
\erem

\brem
The non-proportionality assumption in
Example~\ref{ex:congruence} is essential. If
\[
B=\lambda A
\]
for some $\lambda>0$, then
$\Psi=\lambda\Phi$, and hence
\[
\widehat S_f(\Phi\|\Psi)
=
\lambda f(\lambda^{-1})A.
\]
Thus, the boundary-value formula in
Example~\ref{ex:congruence} reflects the singularity of the two
maps and does not apply in the proportional case.

In particular, let $\eta(t)=t\log t$. If $A$ and $B$ are not
proportional, then
\[
D_{\mathrm{BS}}^{\mathrm{op}}(\Phi\|\Psi)(\omega)
=
+\infty
\]
for every nonzero $\omega\in M_n(\C)_*^+$, because
$\eta'(+\infty)=+\infty$, $\eta(+0)=0$, and $A$ is invertible.
On the other hand, if $B=\lambda A$, then
\[
D_{\mathrm{BS}}^{\mathrm{op}}(\Phi\|\Psi)
=
-(\log\lambda)A.
\]
Consequently, for any two non-proportional matrices
$A,B\in M_n(\C)_{++}$, the corresponding congruence maps have infinite
BS relative entropy, even though their values $A$ and $B$ at the
identity may be arbitrarily close.
\erem

%%%%%%%%%%%%%%%%%%%%%%%%%%%%%%%%%%%%%%%%%%%%%%%%%%%%%%%%%%%%%%%%%%%%%%%%%%%%%%%%%%%%%%%%%%%%%%%%%%%%%%%%%%%%%%%%%%%%%%%%%%%%%%%%%%%%%%%%%%%%%%%%%%%%%%%%%%%%%%%%%%%%%%%%%%%%%%%%%%%%%%%%%%%%%%%%%%%%%%%%%%%%%%%%%%%%%%%%%%%%%%%%%%%%%%%%
%%%%%%%%%%%%%%%%%%%%%%%%%%%%%%%%%%%%%%%%%%%%%%%%%%%%%%%%%%
\section{Finite-index conditional expectations}
%%%%%%%%%%%%%%%%%%%%%%%%%%%%%%%%%%%%%%%%%%%%%%%%%%%%%%%%%%%%%%%%%%%%%%%%%%%%%%%%%%%%%%%%%%%%%%%%%%%%%%%%%%%%%%%%%%%%%%%%%%%%%%%%%%%%%%%%%%%%%%%%%%%%%%%%%%%%%%%%%%%%%%%%%%%%%%%%%%%%%%%%%%%%%%%%%%%%%%%%%%%%%%%%%%%%%%%%%%%%%%%%%%%%%%%%%%%%%%%%%%%%%%%%%%%%%%%%%%%%%%%%%%%%%%%%%%%%%%%%%%%%%%%%%%

\bthm[Finite-index conditional expectations]
\label{thm:finite-index-expectation}
Let $\cN\subseteq\cM$ be an inclusion of $\sigma$-finite factors,
and let
$
E\colon\cM\to\cN
$
be a faithful normal conditional expectation with finite
Jones--Kosaki index
$
\lambda\coloneqq\operatorname{Ind}E<+\infty.
$
Regarding $E$ as an $\cM$-valued CP map via the inclusion
$\cN\subseteq\cM$, for every
$f\in\mathrm{OC}(0,+\infty)$ we have
\[
\widehat S_f(\id_{\cM}\| E)
=
\left\{
\frac{f(\lambda)}{\lambda}
+
\left(1-\frac{1}{\lambda}\right)f(+0)
\right\}1_{\cM}
\]
in $\widehat{\cM}_{\mathrm{lb}}$, with the convention
$
0\cdot(+\infty)=0.
$
\ethm

\bpf
We first note that
\[
\lambda E-\id_{\cM}\in\mathrm{CP}(\cM,\cM),
\qquad\text{equivalently}\qquad
\id_{\cM}\leq_{\mathrm{cp}}\lambda E;
\]
see \cite[Theorem~3.5 and Remark~3.8]{bdh_1988}; see also
\cite[Section~2.5]{popa_2023} for the Pimsner--Popa inequality.
In particular, $\lambda\geq1$.
Therefore, by Proposition~\ref{prop:independent}, we may use $E$ as
the dominating CP map for the pair $(\id_{\cM},E)$.

Choose a faithful normal state 
$\psi\in\cN_*^+$, and set
\[
\varphi\coloneqq\psi\circ E.
\]
Then $\varphi$ is a faithful normal state on $\cM$ satisfying
$
\varphi\circ E=\varphi$.
Consider the GNS representation of $\cM$ associated with $\varphi$,
and write
\[
\cH=L^2(\cM,\varphi),
\qquad
\varphi=\omega_{\xi_0},
\]
where $\xi_0$ is the canonical cyclic and separating vector.
Let $J=J_\varphi$ denote the corresponding modular conjugation.

Define the Jones projection $e_{\cN}\in\cN'$ by
\[
e_{\cN}(x\xi_0)=E(x)\xi_0,
\qquad x\in\cM.
\]
It satisfies
\[
e_{\cN}ae_{\cN}=E(a)e_{\cN},
\qquad a\in\cM.
\]
Set
$
\cM_1\coloneqq\langle\cM,e_{\cN}\rangle$.
Then
$
\cM_1=J\cN'J$.

Let
$
E^{-1}\colon\cN'\to\cM'
$
be Kosaki's operator-valued weight associated with $E$. Since
$\cM$ is a factor, we have
\[
E^{-1}(1)=\lambda I_{\cH}.
\]
The dual conditional expectation
$
E_1\colon\cM_1\to\cM
$
is given by
\[
E_1(z)
=
\lambda^{-1}JE^{-1}(JzJ)J,
\qquad z\in\cM_1,
\]
and satisfies
\[
E_1(xe_{\cN}y)=\lambda^{-1}xy,
\qquad x,y\in\cM.
\]
We refer to \cite[Section~3]{kosaki_1986} for these
basic-construction identities.

Set
\[
\varphi_1\coloneqq\varphi\circ E_1,
\qquad
\cK\coloneqq L^2(\cM_1,\varphi_1).
\]
Since $E_1$ is faithful and normal, $\varphi_1$ is a faithful normal
state on $\cM_1$. Let $\xi_1\in\cK$ be the canonical cyclic and
separating vector implementing $\varphi_1$.
Let $\cM_1$ act on $\cK$ by left multiplication, and denote by $\pi$
the restriction of this representation to $\cM$.

Define an operator $V$ on the dense subspace
$\cM\xi_0\subseteq\cH$ by
\[
V(x\xi_0)
\coloneqq
\lambda^{1/2}e_{\cN}x\xi_1,
\qquad x\in\cM.
\]
For $x,y\in\cM$, we have
\begin{align*}
\langle Vx\xi_0,Vy\xi_0\rangle_{\cK}
&=
\lambda\varphi_1(y^*e_{\cN}x)
=
\lambda\varphi(E_1(y^*e_{\cN}x))\\
&=
\varphi(y^*x)
=
\langle x\xi_0,y\xi_0\rangle_{\cH}.
\end{align*}
Thus $V$ extends uniquely to an isometry from $\cH$ into $\cK$.

Moreover, for $a,x,y\in\cM$, we obtain
\begin{align*}
\langle V^*\pi(a)Vx\xi_0,y\xi_0\rangle_{\cH}
&=
\lambda
\langle ae_{\cN}x\xi_1,e_{\cN}y\xi_1\rangle_{\cK}
=
\lambda
\langle E(a)e_{\cN}x\xi_1,y\xi_1\rangle_{\cK}\\
&=
\lambda\varphi(E_1(y^*E(a)e_{\cN}x))
=
\varphi(y^*E(a)x)\\
&=
\langle E(a)x\xi_0,y\xi_0\rangle_{\cH}.
\end{align*}
Hence
\[
V^*\pi(a)V=E(a),
\qquad a\in\cM.
\]

Since $\operatorname{span}\cM e_{\cN}\cM$ is strongly dense in
$\cM_1$ and $\xi_1$ is cyclic for $\cM_1$, we have
\[
\overline{\pi(\cM)V\cH}
=
\overline{\operatorname{span}}
\{xe_{\cN}y\xi_1\mid x,y\in\cM\}
=
\cK.
\]
Therefore, $(\pi,V,\cK)$ is a minimal Stinespring representation
of $E$.

Set
\[
\cD
\coloneqq
\operatorname{span}
\{xe_{\cN}y\xi_1\mid x,y\in\cM\}\subseteq\cK.
\]
Define a linear map $T_0$ by
\[
T_0
\left(
\sum_{i=1}^n x_i e_{\cN}y_i\xi_1
\right)
\coloneqq
\sum_{i=1}^n x_i y_i e_{\cN}\xi_1.
\]
We show that $T_0$ is well defined and bounded.
Let
\[
\zeta=\sum_{i=1}^n x_i e_{\cN}y_i\xi_1\in\cD,
\qquad
Y=(y_1\xi_0,\ldots,y_n\xi_0)^T\in\cH^n,
\]
and set
\[
G=[x_i^*x_j]_{i,j=1}^n\in M_n(\cM)_+.
\]
By the complete positivity of $\lambda E-\id_{\cM}$ noted at the beginning of the proof,
\[
\lambda^{-1}G\leq E^{(n)}(G).
\]
Then
\begin{align*}
\|T_0\zeta\|_{\cK}^2
&=
\sum_{i,j}
\varphi_1(e_{\cN}y_j^*x_j^*x_i y_i e_{\cN})
=
\lambda^{-1}
\sum_{i,j}
\varphi(y_j^*x_j^*x_i y_i)\\
&=
\lambda^{-1}
\langle GY,Y\rangle_{\cH^n}
\leq
\langle E^{(n)}(G)Y,Y\rangle_{\cH^n}\\
&=
\sum_{i,j}
\varphi(y_j^*E(x_j^*x_i)y_i)
=
\lambda\|\zeta\|_{\cK}^2.
\end{align*}
In particular, if $\zeta=0$, then $T_0\zeta=0$. 
Thus $T_0$ is well defined on $\cD$ and extends uniquely to an operator
$T\in\bB(\cK)$ satisfying
$
\|T\|\leq\lambda^{1/2}$.

For $a,x,y\in\cM$, we have
\[
T\pi(a)(xe_{\cN}y\xi_1)
=
axye_{\cN}\xi_1
=
\pi(a)T(xe_{\cN}y\xi_1).
\]
Since $\cD$ is dense in $\cK$, it follows that
$
T\in\pi(\cM)'$.
Set
\[
A'\coloneqq T^*T\in\pi(\cM)'.
\]
Then
$
0\leq A'\leq\lambda I_{\cK}$.
For $a,x,y\in\cM$, we have
\begin{align*}
\langle V^*A'\pi(a)Vx\xi_0,y\xi_0\rangle_{\cH}
&=
\lambda
\langle T(ae_{\cN}x\xi_1),T(e_{\cN}y\xi_1)\rangle_{\cK}\\
&=
\lambda
\langle axe_{\cN}\xi_1,ye_{\cN}\xi_1\rangle_{\cK}
=
\lambda\varphi_1(e_{\cN}y^*axe_{\cN})\\
&=
\varphi(y^*ax)
=
\langle ax\xi_0,y\xi_0\rangle_{\cH}.
\end{align*}
Thus $A'$ is the RN derivative of $\id_{\cM}$ with respect to
the dominating map $E$.

Since $E_1|_{\cM}=\id_{\cM}$, we have
$
\varphi_1|_{\cM}=\varphi$.
Therefore, the map
\[
x\xi_0\mapsto x\xi_1,
\qquad x\in\cM,
\]
extends uniquely to a unitary from $L^2(\cM,\varphi)$ onto
$\overline{\cM\xi_1}$. We identify these two Hilbert spaces
through this unitary.

Let $P$ be the orthogonal projection from $\cK$ onto
$\overline{\cM\xi_1}$. For $a,x,y\in\cM$, we have
\begin{align*}
\langle P(xe_{\cN}y\xi_1),a\xi_1\rangle_{\cK}
&=
\langle xe_{\cN}y\xi_1,a\xi_1\rangle_{\cK}
=
\varphi_1(a^*xe_{\cN}y)\\
&=
\lambda^{-1}\varphi(a^*xy)
=
\left\langle
\lambda^{-1}xy\xi_1,a\xi_1
\right\rangle_{\cK}.
\end{align*}
Hence
\[
P(xe_{\cN}y\xi_1)
=
\lambda^{-1}xy\xi_1.
\]

We next compute $T^*$. For $a,x,y\in\cM$, we have
\begin{align*}
\langle T^*(ae_{\cN}\xi_1),xe_{\cN}y\xi_1\rangle_{\cK}
&=
\langle ae_{\cN}\xi_1,xye_{\cN}\xi_1\rangle_{\cK}\\
&=
\varphi_1(e_{\cN}y^*x^*ae_{\cN})
=
\lambda^{-1}\varphi(y^*x^*a),
\end{align*}
whereas
\[
\langle a\xi_1,xe_{\cN}y\xi_1\rangle_{\cK}
=
\varphi_1(y^*e_{\cN}x^*a)
=
\lambda^{-1}\varphi(y^*x^*a).
\]
Since $\cD$ is dense in $\cK$, it follows that
\[
T^*(ae_{\cN}\xi_1)=a\xi_1.
\]
Consequently,
\begin{align*}
A'(xe_{\cN}y\xi_1)
&=
T^*T(xe_{\cN}y\xi_1)
=
T^*(xye_{\cN}\xi_1)\\
&=
xy\xi_1
=
\lambda P(xe_{\cN}y\xi_1).
\end{align*}
Since both $A'$ and $\lambda P$ are bounded and this identity
holds on the dense subspace $\cD$, we conclude that
\[
A'=\lambda P.
\]

Furthermore,
\begin{align*}
\langle V^*PVx\xi_0, y\xi_0 \rangle_{\cH}
&=
\lambda\langle Pe_{\cN}x\xi_1, e_{\cN}y\xi_1\rangle_{\cK}
\\
&=
\lambda^{-1}\langle x\xi_1, y\xi_1 \rangle_{\cK}
=
\lambda^{-1}\langle x\xi_0, y\xi_0 \rangle_{\cH}.
\end{align*}
This implies that
\[
V^*PV=\lambda^{-1}I_{\cH}.
\]

The RN derivative of $E$ with respect to itself is $I_{\cK}$.
Therefore, Proposition~\ref{prop:independent} gives
\[
\widehat S_f(\id_{\cM}\| E)
=
V^*\phi_f(A',I_{\cK})V
=
V^*f(A')V.
\]
Since $A'=\lambda P$, the spectral calculus gives
\[
f(A')
=
f(\lambda)P+f(+0)(I_{\cK}-P)
\]
in $\widehat{\bB(\cK)}_{\mathrm{lb}}$. 
Consequently, in the sense of extended lower-semibounded self-adjoint elements,
\begin{align*}
\widehat S_f(\id_{\cM}\| E)
&=
f(\lambda)V^*PV
+
f(+0)V^*(I_{\cK}-P)V\\
&=
\frac{f(\lambda)}{\lambda}I_{\cH}
+
\left(1-\frac{1}{\lambda}\right)f(+0)I_{\cH}.
\end{align*}
\epf

\bcor\label{cor:finite-index-bs}
Under the assumptions of
Theorem~\ref{thm:finite-index-expectation}, let
$f\in\mathrm{OC}(0,+\infty)$ satisfy $f(+0)=0$. 
Then
\[
\widehat S_f(\id_{\cM}\| E)
=
\frac{f(\operatorname{Ind}E)}
{\operatorname{Ind}E}\,1_{\cM}.
\]
In particular,
\[
D_{\mathrm{BS}}^{\mathrm{op}}(\id_{\cM}\| E)
=
(\log\operatorname{Ind}E)1_{\cM}.
\]
Consequently, for every $\varphi\in\cM_*^+$,
\[
D_{\mathrm{BS}}^{\mathrm{op}}(\id_{\cM}\| E)(\varphi)
=
(\log\operatorname{Ind}E)\,\varphi(1_{\cM}).
\]
In particular, for every normal state $\varphi$ on $\cM$,
\[
D_{\mathrm{BS}}^{\mathrm{op}}(\id_{\cM}\| E)(\varphi)
=
\log\operatorname{Ind}E.
\]
\ecor

%%%%%%%%%%%%%%%%%%%%%%%%%%%%%%%%%%%%%%%%%%%%%%%%%%%%%%%%%%%%%%%%%%%%%%%%%%%%%%%%%%%%%%%
%%%%%%%%%%%%%%%%%%%%%%%%%%%%%%%%%%%%%%%%%%%%%%%%%%%%%%%%%%%%%%%%%%%%%%%%%%%%%%%%%%%%%%%
\section{Comparison with the channel divergence of Hollands and Ranallo}
%%%%%%%%%%%%%%%%%%%%%%%%%%%%%%%%%%%%%%%%%%%%%%%%%%%%%%%%%%%%%%%%%%%%%%%%%%%%%%%%%%%%%%%
%%%%%%%%%%%%%%%%%%%%%%%%%%%%%%%%%%%%%%%%%%%%%%%%%%%%%%%%%%%%%%%%%%%%%%%%%%%%%%%%%%%%%%%

Let
$\eta(t)\coloneqq t\log t$.
Let $\cB$ be a unital $*$-algebra,
and let $\varphi,\psi$ be
positive linear functionals on $\cB$.
For $n\in\mathbb N$, denote by $\mathscr X_n(\cB)$ the set of
finite-range step functions
\[
x\colon (1/n,+\infty)\to\cB
\]
satisfying 
$x_t=1_{\cB}$ for $t$ near $1/n$,
and $x_t=0$ for sufficiently large $t$.
Motivated by the variational formula of Hollands and Ranallo
\cite[Proposition~2.17 and Definition~3.4]{h-r_2023},
we use the following formulation:
\begin{align}
D_{\mathrm{BS}}^{\cB}(\varphi\|\psi)
\coloneqq
\sup_{n\in\mathbb N}
\sup_{x\in\mathscr X_n(\cB)}
\Bigg\{
\varphi(1_{\cB})\log n
-
\int_{1/n}^{+\infty}
\left[
\varphi(x_tx_t^*)
+\frac{1}{t}\psi(y_ty_t^*)
\right]
\frac{dt}{t}
\Bigg\},
\label{eq:HR-variational-formula}
\end{align}
where
$y_t\coloneqq1_{\cB}-x_t$.

When $\cB$ is a von Neumann algebra and $\varphi,\psi$ are normal,
the right-hand side of \eqref{eq:HR-variational-formula} is the
variational expression of the BS relative entropy
given in \cite[Proposition~2.17]{h-r_2023}; thus
$D_{\mathrm{BS}}^{\cB}(\varphi\|\psi)$ coincides with the BS relative
entropy used in \cite{h-r_2023}.
For a general unital $*$-algebra $\cB$, suppose that
$\varphi,\psi$ are realized as in Lemma~\ref{lem:hr-pw} below with
$\zeta$ cyclic for $\rho(\cB)$, and let
$\bar\varphi,\bar\psi$ be the normal positive functionals on
$\cR\coloneqq\rho(\cB)''$ given by
$\bar\varphi(y)=\langle Cy\zeta,\zeta\rangle$ and
$\bar\psi(y)=\langle Dy\zeta,\zeta\rangle$.
Since $C,D\in\cR'$ and $\zeta$ is also cyclic for $\cR$,
Lemma~\ref{lem:hr-pw}, applied to both $\cB$ and $\cR$, gives
\[
D_{\mathrm{BS}}^{\cB}(\varphi\|\psi)
=
\phi_\eta(C,D)(\omega_\zeta)
=
D_{\mathrm{BS}}^{\cR}(\bar\varphi\|\bar\psi).
\]
Hence, for $\cB=\cM\odot\cA^{\mathrm{op}}$, the quantity
$D_{\mathrm{BS}}^{\cB}$ appearing in the definition of
$D_{\mathrm{BS}}^{\mathrm{HR}}$ below agrees with the BS relative
entropy of the normal extensions in the sense of
\cite[Proposition~2.17]{h-r_2023}.

\blem
\label{lem:hr-pw}
Let $\cB$ be a unital $*$-algebra.
Suppose that
$\rho\colon\cB\to\bB(\cK)$ is a unital $*$-representation,
$\zeta\in\cK$, and $C,D\in\rho(\cB)'_+$ satisfy
$
C+D=I_{\cK}
$
and 
\[
\varphi(b)
=
\langle C\rho(b)\zeta,\zeta\rangle,
\qquad
\psi(b)
=
\langle D\rho(b)\zeta,\zeta\rangle,
\qquad b\in\cB.
\]
Then
\[
D_{\mathrm{BS}}^{\cB}(\varphi\|\psi)
\leq
\phi_\eta(C,D)(\omega_\zeta).
\]
If $\zeta$ is cyclic for $\rho(\cB)$,
then
\[
D_{\mathrm{BS}}^{\cB}(\varphi\|\psi)
=
\phi_\eta(C,D)(\omega_\zeta).
\]

\elem

\bpf
Set
\[
c\coloneqq
\langle C\zeta,\zeta\rangle
=
\varphi(1_{\cB}).
\]
By \cite[Example~9.8(2), equation~(9.15)]{huw_2022},
we have
\begin{align}
\phi_\eta(C,D)(\omega_\zeta)
=
\sup_{(X,Y)}
\int_0^{+\infty}
\Bigg[
\frac{c}{1+\lambda}
-
\frac{1}{\lambda}
\langle CX(\lambda),X(\lambda)\rangle
-
\langle DY(\lambda),Y(\lambda)\rangle
\Bigg]
\,d\lambda,
\label{eq:huw-bs-variational}
\end{align}
where the supremum is taken over all pairs
$(X(\cdot),Y(\cdot))$ of $\cK$-valued finite-range step functions satisfying
\[
X(\lambda)+Y(\lambda)=\zeta,
\qquad \lambda>0,
\]
such that
$X(\lambda)=0$
for sufficiently small $\lambda>0$,
and
$Y(\lambda)=0$
for sufficiently large $\lambda>0$.
Making the change of variables
$\lambda=t^{-1}$ and setting
\[
u(t)\coloneqq X(t^{-1}),
\qquad
v(t)\coloneqq Y(t^{-1}),
\]
we have
\[
u(t)+v(t)=\zeta,\qquad t>0.
\]
Moreover,
$
u(t)=0
$
for sufficiently large $t>0$,
and
$
v(t)=0
$
for sufficiently small $t>0$.
In particular,
$
u(t)=\zeta
$
for sufficiently small $t>0$.
Then the integral in \eqref{eq:huw-bs-variational} becomes
\begin{align}
\int_0^{+\infty}
\Bigg[
\frac{c}{t(1+t)}
-
\frac{1}{t}
\langle Cu(t),u(t)\rangle
-
\frac{1}{t^2}
\langle Dv(t),v(t)\rangle
\Bigg]
\,dt.
\label{eq:reciprocal-bs-variational}
\end{align}

Choose $n\in\mathbb N$ sufficiently large 
such that
\[
u(t)=\zeta,
\qquad
v(t)=0,
\qquad
0<t\leq1/n.
\]
Splitting the integral in
\eqref{eq:reciprocal-bs-variational} at $1/n$, we obtain
\begin{align}
&
\int_0^{+\infty}
\Bigg[
\frac{c}{t(1+t)}
-
\frac{1}{t}
\langle Cu(t),u(t)\rangle
-
\frac{1}{t^2}
\langle Dv(t),v(t)\rangle
\Bigg]
\,dt
\notag
\\
&=
c\log n
-
\int_{1/n}^{+\infty}
\left[
\langle Cu(t),u(t)\rangle
+
\frac1t
\langle Dv(t),v(t)\rangle
\right]
\frac{dt}{t}.
\label{eq:hr-vector-variational}
\end{align}

For $x\in\mathscr X_n(\cB)$, 
extend $x$ to $(0,+\infty)$ by
setting
\[
x_t=1_{\cB},
\qquad 0<t\leq1/n,
\]
and set $y_t\coloneqq1_{\cB}-x_t$.
If we put
\[
u(t)\coloneqq\rho(x_t^*)\zeta,
\qquad
v(t)\coloneqq\rho(y_t^*)\zeta,
\]
then we have
\[
u(t)+v(t)=\zeta.
\]
Since $C,D\in\rho(\cB)'$, we also have
\[
\langle Cu(t),u(t)\rangle
=
\left\langle
C\rho(x_t^*)\zeta,
\rho(x_t^*)\zeta
\right\rangle
=
\left\langle
C\rho(x_tx_t^*)\zeta,
\zeta
\right\rangle
=
\varphi(x_tx_t^*),
\]
and similarly,
\[
\langle Dv(t),v(t)\rangle
=
\psi(y_ty_t^*).
\]
Consequently, the right-hand side of
\eqref{eq:hr-vector-variational} becomes
\[
\varphi(1_{\cB})\log n
-
\int_{1/n}^{+\infty}
\left[
\varphi(x_tx_t^*)
+
\frac1t\psi(y_ty_t^*)
\right]
\frac{dt}{t},
\]
which is precisely the Hollands--Ranallo variational formula
\eqref{eq:HR-variational-formula}.

Thus, the admissible pairs arising from
$x\in\mathscr X_n(\cB)$ form a subclass of all admissible pairs
in the PW variational formula. Taking the supremum over this
restricted subclass gives
\[
D_{\mathrm{BS}}^{\cB}(\varphi\|\psi)
\leq
\phi_\eta(C,D)(\omega_\zeta).
\]

Conversely, assume that $\zeta$ is cyclic for $\rho(\cB)$.
Let $(u,v)$ be an arbitrary admissible pair in
the PW variational formula.

Choose $0<\delta<R<+\infty$ such that
\[
u(t)=\zeta,
\qquad
v(t)=0,
\qquad
0<t\leq\delta,
\]
and
\[
u(t)=0,
\qquad
v(t)=\zeta,
\qquad
t\geq R.
\]
Take $n\in\mathbb N$ so that
\[
\frac1n<\delta.
\]
Since $u$ is a finite-range step function, after refining its
underlying finite partition, write
\[
u(t)=u_j
\]
on each of finitely many intervals
\[
I_j\subseteq[\delta,R].
\]
By cyclicity, for every $j$ there exist
$b_j^{(k)}\in\cB$ such that
\[
\rho(b_j^{(k)*})\zeta\to u_j.
\]
Define $x^{(k)}\in\mathscr X_n(\cB)$ by setting
$x_t^{(k)}=1_{\cB}$ on $(1/n,\delta)$,
$x_t^{(k)}=b_j^{(k)}$ on $I_j$, and
$x_t^{(k)}=0$ for $t\geq R$.
Then
\[
\rho(x_t^{(k)*})\zeta\to u(t),
\qquad
\rho((1_{\cB}-x_t^{(k)})^*)\zeta\to v(t)
\]
on each interval of the partition.
Since $C$ and $D$ are bounded and the weights $t^{-1}$ and
$t^{-2}$ are integrable on $[\delta,R]$, while the approximating
pairs agree exactly with $(u,v)$ outside this interval, the
corresponding variational expressions converge to the value
associated with $(u,v)$.
Taking the supremum over all
admissible pairs gives
\[
\phi_\eta(C,D)(\omega_\zeta)
\leq
D_{\mathrm{BS}}^{\cB}(\varphi\|\psi).
\]
Together with the opposite inequality proved above, this yields
the desired equality.
\epf

A linear map $\Phi\colon\cM\to\cN$
is called a \emph{channel} if it is unital and CP.

Let $\cA$ be a von Neumann algebra and let
\[
\theta\colon\cN\odot\cA^{\mathrm{op}}
\to\bB(\cH_\theta)
\]
be a binormal unital $*$-representation. Write
\[
L_\theta(x)
\coloneqq\theta(x\otimes1),
\qquad
R_\theta(a^{\mathrm{op}})
\coloneqq
\theta(1\otimes a^{\mathrm{op}}).
\]
Set $\cB\coloneqq\cM\odot\cA^{\mathrm{op}}$.
For a normal channel $\Phi\colon\cM\to\cN$ and a unit vector
$\xi\in\cH_\theta$, define a linear functional on $\cB$ by
\[
\omega_{\Phi,\theta,\xi}
(x\otimes a^{\mathrm{op}})
\coloneqq
\left\langle
L_\theta(\Phi(x))R_\theta(a^{\mathrm{op}})\xi,\xi
\right\rangle,
\qquad
x\in\cM,\ a\in\cA,
\]
and extend it linearly to $\cB$. 
Equivalently,
\[
\omega_{\Phi,\theta,\xi}
=
\omega_\xi\circ\theta\circ
(\Phi\odot\id_{\cA^{\mathrm{op}}}).
\]
The functional $\omega_{\Phi,\theta,\xi}$ is a state on $\cB$.
Indeed, $\omega_{\Phi,\theta,\xi}(1_{\cB})=\|\xi\|^2=1$, and for
$b=\sum_{i=1}^n x_i\otimes a_i^{\mathrm{op}}\in\cB$, since
$L_\theta(\cN)$ and $R_\theta(\cA^{\mathrm{op}})$ commute, we have
\[
\omega_{\Phi,\theta,\xi}(b^*b)
=
\sum_{i,j=1}^n
\bigl\langle
L_\theta\bigl(\Phi(x_i^*x_j)\bigr)\eta_j,\eta_i
\bigr\rangle
\geq0,
\qquad
\eta_j\coloneqq R_\theta(a_j^{\mathrm{op}})\xi,
\]
because $[\Phi(x_i^*x_j)]_{i,j}\in M_n(\cN)_+$ by the complete
positivity of $\Phi$ and $L_\theta$ is a $*$-homomorphism.

Following \cite[Definition~3.4]{h-r_2023}, the Hollands--Ranallo
BS channel divergence is defined by
\[
D_{\mathrm{BS}}^{\mathrm{HR}}(\Phi\|\Psi)
\coloneqq
\sup_{(\cA,\theta,\xi)}
D_{\mathrm{BS}}^{\cB}
\left(
\omega_{\Phi,\theta,\xi}
\middle\|
\omega_{\Psi,\theta,\xi}
\right),
\]
where the supremum is taken over all von Neumann algebras
$\cA$, all binormal unital $*$-representations $\theta$ as above, and
all unit vectors $\xi\in\cH_\theta$.

\begin{theorem}\label{thm:HR-comparison}
Let $\Phi,\Psi\colon\cM\to\cN$ be normal channels and set
\[
m\coloneqq
D_{\mathrm{BS}}^{\mathrm{op}}(\Phi\|\Psi)
\in\widehat{\cN}_{\mathrm{lb}}.
\]
Then, for every von Neumann algebra $\cA$, 
every binormal unital $*$-representation
\[
\theta\colon
\cN\odot\cA^{\mathrm{op}}
\to
\bB(\cH_\theta),
\]
and every unit vector $\xi\in\cH_\theta$, we have
\[
D_{\mathrm{BS}}^{\cB}
\bigl(
\omega_{\Phi,\theta,\xi}
\|
\omega_{\Psi,\theta,\xi}
\bigr)
\leq
\widehat{L_\theta}(m)(\omega_\xi),
\]
where
$\cB\coloneqq\cM\odot\cA^{\mathrm{op}}$.

In particular, if
\[
D_{\mathrm{BS}}^{\mathrm{op}}(\Phi\|\Psi)
\leq c1_{\cN}
\]
for some $c\geq 0$, then
\[
D_{\mathrm{BS}}^{\mathrm{HR}}(\Phi\|\Psi)
\leq c.
\]
\end{theorem}

\bpf
Set
\[
\alpha\coloneqq L_\theta\circ\Phi,
\qquad
\beta\coloneqq L_\theta\circ\Psi.
\]
We regard $\alpha$ and $\beta$ as CP maps from $\cM$ into the
von Neumann algebra $L_\theta(\cN)\subseteq\bB(\cH_\theta)$,
and let
\[
(\pi,V,\cK;A',B')
\]
be a minimal SRN realization of $(\alpha,\beta)$.
Thus,
\[
\alpha(x)=V^*A'\pi(x)V,
\qquad
\beta(x)=V^*B'\pi(x)V,
\qquad
A'+B'=I_{\cK}.
\]

By Lemma~\ref{lem:pi-commutant}, there exists a
unital $*$-representation
$
\pi'\colon
L_\theta(\cN)'
\to
\bB(\cK)
$
such that
\[
\pi'(y')V=Vy',
\qquad y'\in L_\theta(\cN)',
\]
and
\[
\pi'(L_\theta(\cN)')
\subseteq\pi(\cM)',
\qquad
A',B'\in\pi'(L_\theta(\cN)')'.
\]

Since
\[
R_\theta(\cA^{\mathrm{op}})
\subseteq L_\theta(\cN)',
\]
we can define
\[
\widetilde R_\theta(a^{\mathrm{op}})
\coloneqq
\pi'
\bigl(R_\theta(a^{\mathrm{op}})\bigr),
\qquad
a^{\mathrm{op}}\in\cA^{\mathrm{op}}.
\]
Hence, the commuting representations $\pi$ and
$\widetilde R_\theta$ induce a unital $*$-representation
\[
\rho\colon
\cM\odot\cA^{\mathrm{op}}
\to\bB(\cK)
\]
given by
\[
\rho(x\otimes a^{\mathrm{op}})
\coloneqq
\pi(x)\widetilde R_\theta(a^{\mathrm{op}}).
\]
Since $A'$ and $B'$ commute with both $\pi(\cM)$ and
$\widetilde R_\theta(\cA^{\mathrm{op}})$, we have
\[
A',B'\in\rho(\cB)'.
\]

Put
\[
\zeta\coloneqq V\xi\in\cK.
\]
For $x\in\cM$ and
$a^{\mathrm{op}}\in\cA^{\mathrm{op}}$, we have
\begin{align*}
\left\langle
A'\rho(x\otimes a^{\mathrm{op}})\zeta,
\zeta
\right\rangle
&=
\left\langle
A'\pi(x)\widetilde{R}_\theta(a^{\mathrm{op}})V\xi,
V\xi
\right\rangle
=
\left\langle
A'\pi(x)\pi'(R_\theta(a^{\mathrm{op}}))V\xi,
V\xi
\right\rangle
\\
&=
\left\langle
V^*A'\pi(x)VR_\theta(a^{\mathrm{op}})\xi,
\xi
\right\rangle
=
\left\langle
L_\theta(\Phi(x))
R_\theta(a^{\mathrm{op}})\xi,
\xi
\right\rangle
\\
&=
\omega_{\Phi,\theta,\xi}
(x\otimes a^{\mathrm{op}}).
\end{align*}
Similarly,
\[
\left\langle
B'\rho(x\otimes a^{\mathrm{op}})\zeta,
\zeta
\right\rangle
=
\omega_{\Psi,\theta,\xi}
(x\otimes a^{\mathrm{op}}).
\]

By Lemma~\ref{lem:hr-pw} and 
the postcomposition monotonicity in Theorem~\ref{thm:postcomposition}, we obtain
\begin{align*}
D_{\mathrm{BS}}^{\cB}
\bigl(
\omega_{\Phi,\theta,\xi}
\|
\omega_{\Psi,\theta,\xi}
\bigr)
&\leq
\phi_\eta(A',B')(\omega_\zeta)
\\
&=
\widehat S_\eta
(L_\theta\circ\Phi\|L_\theta\circ\Psi)
(\omega_\xi)
\leq
\widehat{L_\theta}(m)(\omega_\xi).
\end{align*}
 
Suppose now that
\[
m\leq c1_{\cN}.
\]
Since $\widehat{L_\theta}$ is order-preserving and
$L_\theta$ is unital,
\[
\widehat{L_\theta}(m)
\leq
cI_{\cH_\theta}.
\]
Hence, for every admissible triple
$(\cA,\theta,\xi)$, we have
\[
D_{\mathrm{BS}}^{\cB}
\bigl(
\omega_{\Phi,\theta,\xi}
\|
\omega_{\Psi,\theta,\xi}
\bigr)
\leq c.
\]
Taking the supremum over all admissible triples gives
\[
D_{\mathrm{BS}}^{\mathrm{HR}}(\Phi\|\Psi)
\leq c.
\]
\epf

For $m\in\widehat{\cN}_+$, define its extended norm by
\[
\|m\|_{\mathrm{ext}}
\coloneqq
\sup\{
m(\varphi)
\mid
\varphi\in\cN_*^+,\ \varphi(1_{\cN})=1
\}
\in[0,+\infty].
\]
Equivalently,
\[
\|m\|_{\mathrm{ext}}
=
\inf\{
c\geq0
\mid
m\leq c1_{\cN}
\},
\]
where the infimum of the empty set is understood as $+\infty$.
Indeed, by the definition of the order on $\widehat{\cN}_+$,
\[
m\leq c1_{\cN}
\quad\Longleftrightarrow\quad
m(\varphi)\leq c\,\varphi(1_{\cN})
\qquad
\varphi\in\cN_*^+.
\]
Normalizing every nonzero $\varphi\in\cN_*^+$ gives the stated
equivalence. Moreover, if the supremum is finite, then
$m\leq c1_{\cN}$ for some $c<+\infty$; the spectral description in
Proposition~\ref{prop:resolution} then shows that $m$ has no infinite
part and is represented by a bounded positive element of $\cN$.
In particular,
\[
\|m\|_{\mathrm{ext}}
=
\begin{cases}
\|m\| & (m\in\cN_+),\\
+\infty & (m\in\widehat{\cN}_+\setminus\cN_+).
\end{cases}
\]

\bthm\label{thm:hr-norm}
Let $\cM$ be a von Neumann algebra, let $\cN$ be a $\sigma$-finite
von Neumann algebra, and let
$\Phi,\Psi\colon\cM\to\cN$ be normal channels. Then
\[
D_{\mathrm{BS}}^{\mathrm{HR}}(\Phi\|\Psi)
=
\left\|
D_{\mathrm{BS}}^{\mathrm{op}}(\Phi\|\Psi)
\right\|_{\mathrm{ext}}.
\]
\ethm

\bpf
Set
\[
m\coloneqq
D_{\mathrm{BS}}^{\mathrm{op}}(\Phi\|\Psi)
\in\widehat{\cN}_{\mathrm{lb}}.
\]
Let
\[
(\cN,L^2(\cN),J,\cP)
\]
be a standard form
and let
\[
(\pi,V,\cK;A',B')
\]
be a minimal SRN realization of $(\Phi,\Psi)$ with respect to this standard
representation.

We first show that $m$ is positive.
Since
\[
\eta(t)\geq t-1, \qquad t>0,
\]
and since $\eta(+0)=0$ and
$\eta'(+\infty)=+\infty$, we have
\[
\phi_\eta(s,t)\geq s-t,
\qquad s,t\geq0.
\]
Hence, by the order-preserving property of the PW
functional calculus,
\[
\phi_\eta(A',B')\geq A'-B'.
\]
Therefore,
\[
m
=
V^*\phi_\eta(A',B')V
\geq
V^*(A'-B')V
=
\Phi(1_{\cM})-\Psi(1_{\cM})
=
0,
\]
because $\Phi$ and $\Psi$ are unital. Thus,
$m\in\widehat{\cN}_+$.

Let $(\cA,\theta,\xi)$ be an admissible triple in the definition of
$D_{\mathrm{BS}}^{\mathrm{HR}}(\Phi\|\Psi)$. 
Since $L_\theta$ is unital and normal, and $\xi$ is a unit vector,
$\omega_\xi\circ L_\theta$ is a normal state on $\cN$.
By
Theorem~\ref{thm:HR-comparison},
\[
D_{\mathrm{BS}}^{\cB}
\bigl(
\omega_{\Phi,\theta,\xi}
\|
\omega_{\Psi,\theta,\xi}
\bigr)
\leq
\widehat{L_\theta}(m)(\omega_\xi)
=
m(\omega_\xi\circ L_\theta)
\leq
\|m\|_{\mathrm{ext}},
\]
where $\cB=\cM\odot\cA^{\mathrm{op}}$.
Therefore,
\[
D_{\mathrm{BS}}^{\mathrm{HR}}(\Phi\|\Psi)\leq\|m\|_{\mathrm{ext}}.
\]

We next prove the reverse inequality. 
Let
$\psi\in\cN_*^+$ be a faithful normal state, and let
$\xi_\psi\in\cP$ be the unit vector implementing $\psi$. Consider
the standard binormal $*$-representation
\[
\lambda\colon
\cN\odot\cN^{\mathrm{op}}
\to
\bB(L^2(\cN))
\]
defined by
\[
\lambda(a\otimes b^{\mathrm{op}})
=
aJb^*J,
\qquad a,b\in\cN.
\]
Then we have
\[
L_\lambda(a)=a,
\qquad
R_\lambda(b^{\mathrm{op}})=Jb^*J.
\]

By Lemma~\ref{lem:pi-commutant}, there exists a
unital $*$-representation
\[
\pi'\colon\cN'\to\bB(\cK)
\]
such that
\[
\pi'(y')V=Vy',
\qquad y'\in\cN',
\]
and
\[
\pi'(\cN')\subseteq\pi(\cM)',
\qquad
A',B'\in\pi'(\cN')'.
\]

Set $\cB\coloneqq\cM\odot\cN^{\mathrm{op}}$.
Since $Jb^*J\in\cN'$ for $b\in\cN$ and
$\pi'(\cN')\subseteq\pi(\cM)'$,
we may define a unital $*$-representation
\[
\rho
\colon
\cM\odot\cN^{\mathrm{op}}
\to
\bB(\cK)
\]
by
\[
\rho(a\otimes b^{\mathrm{op}})
\coloneqq
\pi(a)\pi'(Jb^*J),
\qquad
a\in\cM,\quad b\in\cN.
\]
Since $A', B'\in\pi(\cM)'\cap\pi'(\cN')'$,
we have 
\[
A', B'\in\rho(\cB)'.
\]
Put
\[
\zeta_\psi\coloneqq V\xi_\psi\in\cK.
\]

For $a\in\cM$ and $b\in\cN$, we have
\begin{align*}
\left\langle
A'\rho(a\otimes b^{\mathrm{op}})\zeta_\psi,
\zeta_\psi
\right\rangle
&=
\left\langle
V^*A'\pi(a)\pi'(Jb^*J)V\xi_\psi,
\xi_\psi
\right\rangle
\\
&=
\left\langle
\Phi(a)Jb^*J\xi_\psi,
\xi_\psi
\right\rangle
=
\omega_{\Phi,\lambda,\xi_\psi}
(a\otimes b^{\mathrm{op}}).
\end{align*}
Similarly,
\[
\left\langle
B'\rho(a\otimes b^{\mathrm{op}})\zeta_\psi,
\zeta_\psi
\right\rangle
=
\omega_{\Psi,\lambda,\xi_\psi}
(a\otimes b^{\mathrm{op}}).
\]
Since $\psi$ is faithful, $\xi_\psi$ is cyclic for $\cN'$.
Hence,
\begin{align*}
\overline{\rho(\cB)\zeta_\psi}
&=
\overline{
\operatorname{span}}
\{
\pi(x)\pi'(y')V\xi_\psi
\mid
x\in\cM,\ y'\in\cN'
\}
\\
&=
\overline{
\operatorname{span}}
\{
\pi(x)V(y'\xi_\psi)
\mid
x\in\cM,\ y'\in\cN'
\}
\\
&=
\overline{\pi(\cM)VL^2(\cN)}
=
\cK.
\end{align*}
Thus $\zeta_\psi$ is cyclic for $\rho(\cB)$. Therefore,
by Lemma~\ref{lem:hr-pw},
we obtain
\[
D_{\mathrm{BS}}^{\cB}
\bigl(
\omega_{\Phi,\lambda,\xi_\psi}
\|
\omega_{\Psi,\lambda,\xi_\psi}
\bigr)
=
\phi_\eta(A',B')(\omega_{\zeta_\psi})
=
m(\omega_{\xi_\psi})=m(\psi).
\]

The triple
$(\cN,\lambda,\xi_\psi)$ 
is admissible in the definition of the
Hollands--Ranallo divergence. Hence,
\[
D_{\mathrm{BS}}^{\mathrm{HR}}(\Phi\|\Psi)
\geq
m(\psi)
\]
for every faithful normal state $\psi$ on $\cN$.

Fix a faithful normal state $\psi_0$ on $\cN$.
Let $\psi$ be an arbitrary normal state on $\cN$ and, for
$0<\varepsilon<1$, set
\[
 \psi_\varepsilon
 \coloneqq
 (1-\varepsilon)\psi+\varepsilon\psi_0.
\]
Then $\psi_\varepsilon$ is faithful
and $\psi_\varepsilon\to\psi$ in norm.
By the lower semicontinuity of
$m\in\widehat{\cN}_+$ and the preceding argument, we have
\[
 m(\psi)
 \leq
 \liminf_{\varepsilon\downarrow0}m(\psi_\varepsilon)
 \leq
 D_{\mathrm{BS}}^{\mathrm{HR}}(\Phi\|\Psi).
\]
Therefore,
\[
D_{\mathrm{BS}}^{\mathrm{HR}}(\Phi\|\Psi)
\geq
\sup\{
m(\omega)
\mid
\omega\in\cN_*^+,
\ \omega(1_{\cN})=1
\}
=
\|m\|_{\mathrm{ext}}.
\]
\epf

For an inclusion $\cN\subseteq\cM$ of factors with finite index, there exists a unique
conditional expectation $E_0\colon\cM\to\cN$ minimizing $\operatorname{Ind}E$, called the
minimal conditional expectation \cite{hiai_1988}; we write
$[\cM:\cN]\coloneqq\operatorname{Ind}E_0$ for the minimal index.

\bcor[cf.~{\cite[Proposition~3.21]{h-r_2023}}]
Let $\cN\subseteq\cM$ be an inclusion of $\sigma$-finite factors
with finite index, and let
$E_0\colon\cM\to\cN$ be the minimal conditional expectation.
Regarding $E_0$ as an $\cM$-valued channel through the inclusion
$\cN\subseteq\cM$, we have
\[
D_{\mathrm{BS}}^{\mathrm{HR}}
(\id_{\cM}\|E_0)
=
\log[\cM:\cN].
\]
\ecor

\bpf
By Corollary~\ref{cor:finite-index-bs} applied to $E=E_0$,
\[
D_{\mathrm{BS}}^{\mathrm{op}}(\id_{\cM}\| E_0)
=
\log[\cM:\cN]\,1_{\cM}.
\]
Hence Theorem~\ref{thm:hr-norm} gives
\[
D_{\mathrm{BS}}^{\mathrm{HR}}(\id_{\cM}\|E_0)
=
\left\|
\log[\cM:\cN]\,1_{\cM}
\right\|_{\mathrm{ext}}
=
\left\|
\log[\cM:\cN]\,1_{\cM}
\right\|
=
\log[\cM:\cN].
\]
\epf

\subsection*{AI statement}

The author used ChatGPT (OpenAI)
and Claude (Anthropic) for English-language editing, to explore possible
proof strategies, and to check the logical structure of the arguments.
The author has verified all proofs and takes full
responsibility for the content of this paper.

%%%%%%%%%%%%%%%%%%%%%%%%%%%%%%%%%%%%%%%%%%%%%%%%%%%%%%%%%%%%%%%%%%%%%%%%%%%%%%%%%%%%%%%%%%%%%%%%%%%%%%%%%%%%%%%%%%%%%%%%%%%%%%%%%%%%%%%%%%%%%%%%%%%%%%%%%%%%%%%%%%%%%%%%%%%%%%%%%%%%%%%%%%%%%%%%%%%%%%%%%%%%%%%%%%%%%%%%%%%%%%%%%%%%%%%%%%%%%

\end{document}